\documentclass[12pt]{amsart}
\usepackage{amsmath}
\usepackage{amsthm}
\usepackage{latexsym}
\usepackage{amssymb}
\DeclareMathOperator{\Char}{Char}

\DeclareMathOperator{\re}{Re}
\DeclareMathOperator{\im}{Im}

\DeclareMathOperator{\dist}{dist}
\DeclareMathOperator{\hess}{Hess}

\DeclareMathOperator{\codim}{codim}
\DeclareMathOperator{\spec}{Spec}
\DeclareMathOperator{\tr}{Tr}

\DeclareMathOperator{\id}{Id}
\DeclareMathOperator{\vol}{Vol}
\DeclareMathOperator{\diam}{diam}
\newtheorem{theorem}{Theorem}
\newtheorem{proposition}{Proposition}
\newtheorem{lemma}{Lemma}

\newtheorem{example}{Example}
\newcounter{obsctr}

\newtheorem{remark}{Remark}
\renewcommand{\thetheorem}{\thesection.\arabic{theorem}}
\renewcommand{\theproposition}{\thesection.\arabic{proposition}}
\renewcommand{\thelemma}{\thesection.\arabic{lemma}}
\renewcommand{\thedefinition}{\thesection.\arabic{definition}}
\renewcommand{\thecorollary}{\thesection.\arabic{corollary}}
\renewcommand{\theequation}{\thesection.\arabic{equation}}
\begin{document}
\baselineskip 18pt
\def\R {{\mathbb{R}}}
\def\N {{\mathbb{N}}}
\def\C {{\mathbb{C}}}
\def\Z {{\mathbb{Z}}}
\def\phi{\varphi}
\def\epsilon{\varepsilon}
%
\def\tb#1{\|\kern -1.2pt | #1 \|\kern -1.2pt |} 
\def\Qed{\qed\par\medskip\noindent}
%
\title{Analytic Hypoellipticity in the Presence \\[5pt]
of Lower Order Terms}     

\author{Paolo Albano, Antonio Bove}
\address{Dipartimento di Matematica, 
Universit\`a di Bologna, Piazza
di Porta San Donato 5, 40127 Bologna, Italy}
\email{albano@dm.unibo.it}
\address{Dipartimento di Matematica, 
Universit\`a di Bologna, Piazza
di Porta San Donato 5, 40127 Bologna, Italy}
\email{Antonio.Bove@bo.infn.it} 
\author{David S. Tartakoff}
\address{Department of Mathematics, University
of Illinois at Chicago, m/c 249, 851 S.
Morgan St., Chicago IL  60607, USA}
\email{dst@uic.edu}
\date{\today}
\begin{abstract}
  We consider a second order operator with analytic coefficients whose
  principal symbol vanishes exactly to order two on a symplectic real
  analytic manifold. We assume that the first (non degenerate)
  eigenvalue vanishes on a symplectic submanifold of the
  characteristic manifold. In the $ C^{\infty} $ framework this
  situation would mean a loss of 3/2 derivatives (see
  \cite{Helffer1}). We prove that this operator is analytic
  hypoelliptic. 
\par  
The main tool is the FBI transform. A case in which 
$C^\infty$ hypoellipticity fails is also discussed.   
\end{abstract}
\maketitle
\pagestyle{myheadings}
\markboth{P. Albano, A. Bove and D. S. Tartakoff}{Analytic
  hypoellipticity and lower order terms}
\tableofcontents

\section{Introduction}
\renewcommand{\theequation}{\thesection.\arabic{equation}}
\setcounter{equation}{0}
\setcounter{theorem}{0}
\setcounter{proposition}{0}  
\setcounter{lemma}{0}
\setcounter{corollary}{0} 
\setcounter{definition}{0}

It is well known that the hypoellipticity of a partial differential
operator heavily depends on the lower order terms, both in the $
C^{\infty} $ and in the analytic category, as it is shown, e.g. in the
$ C^{\infty} $ category, by the papers \cite{Sj-74} and \cite{BGH}.

The study of $ C^{\infty} $ hypoellipticity has produced a number of
results characterizing the loss of derivatives. We say that the
hypoelliptic operator $ P $ of order $ m $ loses $ q $ derivatives if
whenever $ Pu \in H^{s} $ we have that $ u \in H^{s+m-q} $. 

In this paper we are concerned with the analytic hypoellipticity of a
class of second order operators losing more than one derivative.
The minimal loss for our class, basically the
Boutet-Grigis-Helffer class, is one, provided certain conditions
on the lower order terms are satisfied. 

If the lower order terms do not satisfy the hypoellipticity
conditions, which means they take values in certain discrete sets,
then the operator may fail to be hypoelliptic. A celebrated example is
the boundary Kohn Laplacian $ \Box_{b} $ on the Heisenberg group on
functions. E.M.~Stein, \cite{S}, has shown that if we add a non zero
complex constant to it then we obtain a hypoelliptic operator which is
also analytic hypoelliptic. Later Kwon, using the concatenation method
due to L.~Boutet de Monvel and F.~Treves, has generalized Stein's
result. See also \cite{Tar04} for a purely $L^2$ proof of
Stein's result.

At least in a formal way we may describe our setting as
follows. Consider a second order operator with double
characteristics. To keep it simple we may assume that its principal
part, homogeneous of degree two, is positive semi-definite, which
occurs e.g. in the case of a sum of squares of real vector
fields. Further assume that the characteristic set is a symplectic
manifold on which the principal part, $ p_{2} $, vanishes exactly
to second order. This means that the kernel of the Hessian matrix
of $ p_{2}(x, \xi) $ with respect to $ (x, \xi) $ is exactly the
tangent space to the characteristic manifold at $ (x, \xi) $. We know
that $ p_{2}(x, D) $ has then an unbounded discrete spectrum and we
may think of the $ C^{\infty} $-hypoellipticity conditions of
Boutet-Treves, H\"ormander and Treves as conditions prescribing that
the lower order terms do not add up to the elements of the spectrum to
hit zero. 

In the symplectic case, the spectrum of $ p_{2} $ is given by a set of
functions (symbols) defined on the characteristic manifold, which can
be thought of as a cotangent space in a smaller dimension. Thus we may
think of the eigenvalues of $ p_{2} $, and also of the eigenvalues of
the whole operator $ P $, as symbols defined on the characteristic
manifold. We have $ C^{\infty} $ hypoellipticity when these symbols
are elliptic, but also in other cases.

We can roughly state Kwon's result by saying that if the
principal symbol (of order one) of an eigenvalue
vanishes identically on the characteristic set, but nethertheless is
elliptic of an arbitrary order less than one, then $ P $ is
hypoelliptic.

In \cite{Helffer1} B. Helffer has shown that $ C^{\infty} $
hypoellipticity holds if, roughly, one of the eigenvalues
of the operator degenerates on a submanifold of the characteristic
manifold and is actually a complex symbol of principal type such that
its Poisson bracket with its complex conjugate does not vanish. In
this case he proved that there is a loss of 3/2 derivatives. For the
proof he constructed a parametrix of the operator following Sj\"ostrand
and Boutet de Monvel.

The purpose of the present paper is to study the analytic
hypoellipticity of such an operator. For this we assume that both the
characteristic manifold of the principal part and the characteristic
manifold of the eigenvalue are symplectic and that the
restriction of the symplectic 2-form to these manifolds has constant
rank. 

For the sake of simplicity we consider the lowest eigenvalue; it
degenerates on a symplectic submanifold of the characteristic manifold
and its principal symbol satisfies the same Poisson bracket
condition needed for the loss of 3/2 derivatives, i.e. it is a complex
principal type operator. 

We stress the fact that it is important that both manifolds involved
are symplectic, otherwise we cannot hope for better than Gevrey 4
hypoellipticity, due to a well known connection between the symplectic
character of the characteristic set and analytic hypoellipticity
(\cite{Treves}).

The method we use is inspired by the work of J.~Sj\"ostrand
\cite{Sj-83}. The main idea is to deduce suitable estimates for the
FBI transform of the 
solution $ u $ of the equation $ Pu = f $, when $ f $ is analytic and
then use this to to obtain the $ WF_{a} $ of $ u. $ In order to
deduce {\it a priori} estimates, we first prove certain a priori estimates for the
localized operator and then lift these ``tangential'' estimates to
(micro)local estimates for $ P $. The deduction of the a priori
estimates for the localized operator is done by constructing an
approximate parametrix. The conclusion follows using a
deformation argument for the weight function related to the FBI
phase. Sj\"ostrand has proved that this argument holds under certain
low regularity assumptions which we can actually do without.

\section{Statement of the result}\label{s:sr}
\renewcommand{\theequation}{\thesection.\arabic{equation}}
\setcounter{equation}{0}
\setcounter{theorem}{0}
\setcounter{proposition}{0}  
\setcounter{lemma}{0}
\setcounter{corollary}{0} 
\setcounter{definition}{0}

Let
$$ 
P (x, D) = p_{2}(x, D) + p_{1}(x, D) + \cdots
$$
be a (properly supported) real analytic second order pseudodifferential
operator. We assume 

\begin{itemize}
\item[$(H1)$]
The principal symbol $ p_{2}(x, \xi)  $ of $ P $ vanishes exactly of
order 2 on a symplectic real analytic submanifold $ \Sigma_{1} \subset
T^{*}\R^{n}\setminus\{0\} $. 
\end{itemize}
Another way of stating the above assumption amounts to saying that
$ p_{2} (x, \xi) \sim
d_{\Sigma_{1}}^{2}(x, \xi) |\xi|^{2} $, where $ d_{\Sigma_{1}} $ is
the (homogeneous of degree 0) distance to $ \Sigma_{1} $.
\begin{itemize}
\item[$(H2)$]
The localized principal symbol takes values in a proper cone $ \Gamma
$ of the complex plane $\C$. 
\end{itemize}

In the codimension two case $(H2)$ is equivalent to saying that the
winding number of the localized principal symbol is zero.

We recall the definition of the sub-principal symbol and the
positive trace of the operator $P$:   
$$ 
p^{s} (x, \xi) = p_{1}(x, \xi) + \frac{i}{2} \sum_{j=1}^{n}
\frac{\partial^{2} p_{2}(x, \xi) }{\partial x_{j} \partial \xi_{j}}
$$
and 
$$ 
\tr^{+}p_{2}(\rho) = \sum_{\substack{\mu \in \Gamma \\ \pm \mu \in
    \spec(F(\rho))}} \mu(\rho), 
$$
$ F(\rho) $, $ \rho \in \Sigma_{1},$ being the Hamiltonian matrix of $
p_{2} $, i.e.  
$$ 
\sigma(v, F(\rho) w) = \frac{1}{2} \langle \hess p_{2}(\rho) v, w\rangle,
\qquad \forall \rho \in \Sigma_{1}\, , \quad \forall  v, w \in
T_{\rho}T^{*}\R^{n}
\,, 
$$
where $\sigma$ denotes the symplectic form; for example with respect
to the canonical coordinates $(x,\xi)$, $\sigma = d\xi \wedge dx$.
Set 
$$ 
q(\rho) = p^{s}(\rho) + \tr^{+} p_{2}(\rho), \qquad \rho \in
\Sigma_{1}\, .
$$
Now, we are ready to state our third assumption. 
\begin{itemize}

\item[$(H3)$]
There exists a symplectic real analytic submanifold $ \Sigma_{2}
\subset \Sigma_{1} $, with $ \codim_{\Sigma_{1}} \Sigma_{2} = 2 $, such
that $\Char(q) = \Sigma_{2}$ and 
\end{itemize}
\begin{equation}\label{eq:helf} 
\frac{1}{i} \{ \bar{q}\, ,\,  q \} >
0\quad \text{on}\quad \Sigma_{2}\, .
\end{equation} 
Our result can be formulated as  
\begin{theorem}\label{t:1}  
Under the above assumptions $(H1)$, $(H2)$ and  $(H3)$ the
operator $ P $ is analytic hypoelliptic.
\end{theorem}

\begin{example}
Denote by $ x = (t, y, s) \in \R^{\nu} \times \R \times \R^{n - \nu
  -1} $ and consider the operator  
$$ 
P  = \sum_{j=1}^{\nu} D_{t_{j}}^{2} + |t|^{2} |D_{s}|^{2} - \nu |D_{s}|
+ D_{y} + i y |D_{s}|.
$$
Then $ \Sigma_{1} = \{t = \tau = 0\} $ and 
$ \Sigma_{2} = \{ t = \tau = 0, y = \eta = 0, \sigma \neq 0\} $. 
It is also easy to check Assumptions $(H1)$--$(H3)$. 
The above theorem then implies that $ P $ is
analytic hypoelliptic.
\end{example}

\section{The localized operator}\label{s:3}
\renewcommand{\theequation}{\thesection.\arabic{equation}}
\setcounter{equation}{0}
\setcounter{theorem}{0}
\setcounter{proposition}{0}  
\setcounter{lemma}{0}
\setcounter{corollary}{0} 
\setcounter{definition}{0}

The purpose of this section is to deduce suitable \textit{a priori}
estimates for a localized operator. More precisely, let $ x \in
\R^{ \nu +1} $, $ x = (t, y) $, $ t \in \R^{\nu } $, $ y \in \R $ and denote
by $ \Sigma_{1} = \{(0, y; 0, \eta)| \eta \neq 0\} $; consider the
operator
\begin{equation}
\label{eq:2.1}
\tilde{P}(x, D_{x}) = \tilde{p}_{2}(t, D_{t}) + \lambda^{-1/2}
\tilde{p}_{1}(y, D_{y}),
\end{equation}
where $ \lambda > 1 $ is a large parameter. Here $ \tilde{p}_{2} $ is
defined by
$$ 
\tilde{p}_{2}(t, \tau) = \sum_{|\alpha+\beta|=2}
a_{\alpha\beta}t^{\alpha} \tau^{\beta} - \mu_{0},
$$
for $ a_{\alpha\beta} \in \C $, where $ \mu_{0} $ is a complex
constant defined by
\begin{equation}
\label{eq:2.3}
\mu_{0} = \tr^{+}\tilde{p}_{2} + \frac{i}{2}
\sum_{j=1}^{n} \frac{\partial^{2} \tilde{p}_{2}}{\partial t_{j}
  \partial \tau_{j}}.
\end{equation}
We point out that $ \mu_{0} $ is invariantly defined on $ \Sigma_{1}
$.

Furthermore $ \tilde{p}_{1} $ in (\ref{eq:2.1}) is a linear form in
the variables $ (y, \eta) $:
\begin{equation}
\label{eq:2.3bis}
\tilde{p}_{1}(y, \eta) = \tilde{\ell}_{1}(y,\eta).
\end{equation}
As a consequence of our assumptions on the non-localized operator we
have that $ \tilde{P} $ satisfies the following requests:
\begin{itemize}
\item[$(a)$]
$ \tilde{p}_{2} + \mu_{0} $, as a quadratic form, has a numerical
range which is a sector, $ \Gamma $, in $ \C $ of amplitude strictly
less than $ \pi 
$. 
\item[$(b)$]
$ \tilde{p}_{2} + \mu_{0} $ is a non degenerate complex quadratic
form, i.e. $ \tilde{p}_{2}(t, \tau) + \mu_{0} = 0 $ implies $ t = \tau = 0 $.
\item[$(c)$]
Condition (\ref{eq:helf}) implies that
\begin{equation}
\label{eq:2.2}
\frac{1}{i} \{ \bar{\tilde{p}}_{1}(y, \eta), \tilde{p}_{1}(y, \eta)
\}_{\big|_{\Sigma_{2}}} > 0.
\end{equation}
\end{itemize}

\subsection{The parametrix for $\tilde{P}$}
\label{ss:1} 
\renewcommand{\theequation}{\thesubsection.\arabic{equation}}
\setcounter{equation}{0}
\setcounter{theorem}{0}
\setcounter{proposition}{0}  
\setcounter{lemma}{0}
\setcounter{corollary}{0} 
\setcounter{definition}{0}
%
\renewcommand{\thetheorem}{\thesubsection.\arabic{theorem}}
\renewcommand{\theproposition}{\thesubsection.\arabic{proposition}}
\renewcommand{\thelemma}{\thesubsection.\arabic{lemma}}
\renewcommand{\thedefinition}{\thesubsection.\arabic{definition}}
\renewcommand{\thecorollary}{\thesubsection.\arabic{corollary}}
\renewcommand{\theequation}{\thesubsection.\arabic{equation}}
\renewcommand{\theremark}{\thesubsection.\arabic{remark}}

In this section we construct an approximate parametrix for 
$ \tilde{P} $  in (\ref{eq:2.1}). 
We basically follow Sj\"ostrand's ideas (\cite{Sj-74}, see also
\cite{Helffer1}).

Let $ e_{0*}(t) $ denote the null eigenfunction of $
\tilde{p}_{2}^{*}(t, D_{t}) $, the formal $ L^{2} $ adjoint of $
\tilde{p}_{2}(t, D_{t}) $; $ e_{0*}(t) $ is a rapidly decreasing
function defined in $ \R^{ \nu } $, which we suppose normalized, i.e. $ \|
e_{0*} \|_{L^{2}(\R^{\nu})} = 1$. Furthermore $ e_{0*} $ can be extended
to an entire function in $ \C^{ \nu } $. Similarly, we denote by
$e_{0}(t)$ the normalized null eigenfunction of $\tilde{p}_{2}(t,
D_{t}) $. 

We remark that the eigenfunctions $ e_{0} $ and $ e_{0*} $ have
the same parity with respect to $ t \in \R^{\nu} $. Hence $ \langle
e_{0*}, e_{0} \rangle_{\R_{t}^{\nu}} \neq 0 $.

We need the following operators \medskip
$$
\begin{array}{ll}
 R^-:L^2 (\R_y)\to L^2 (\R^{\nu +1}_{(t,y)})\, , \qquad   
&  R^-u (t,y) =e_{0}(t) u(y)\, , 
\\[16pt]
 R_*^+: L^2 (\R^{\nu+1}_{(t,y)})\to L^2 (\R_y)\, , \qquad  
&R_*^+u(y) =\int_{\R^\nu } \overline{e_{0*}(t)} u(t,y)\, dt.
\end{array}
$$
\medskip

\noindent
We observe that the operator $R^- R_*^+$ can 
be also realized as a pseudodifferential operator of Hermite type with 
respect to the $t$-variables.
We use Weyl quantized pseudodifferential operators. One can show that  
\begin{multline*}
(R^- R_*^+) u(t,y)=\tilde{h}^w (t,D_t) u(t,y)
\\
=\frac 1 {(2\pi )^\nu} \int \int 
e^{i \langle t-s,\tau\rangle }  \tilde{h}\left ( \frac {t+s} 2 ,\tau \right ) u(s,y)
\, ds \, d\tau 
\end{multline*}
where $\tilde{h}$ is a rapidly decreasing symbol in the
variables $t$ and $\tau$ and is defined as 
$$
\tilde{h}(t,\tau ) = 2^\nu \int e^{-2i \langle \tau , s \rangle } e_0 (t+s) 
\overline{e_{0*}(t-s)}\, ds\, .
$$
Furthermore there exists an operator, in the $ t $
variables,  
$$
\tilde{F}: L^2 (\R^{\nu+1}_{(t,y )})\to L^2 (\R^{\nu+1}_{(t,y )})\, , 
$$ 
such that, still denoting by $\tilde{F}$ the Weyl symbol of $\tilde{F}$,  
\begin{equation}\label{eq:f}
\tilde{F} \# \tilde{p}_2 = 1 - \frac {\tilde{h}}{c_0}\, ,
\end{equation} 
where 
$$
c_0=R_*^+R^- \, .
$$
Here $\#$ denotes the Weyl composition.
Moreover, we have
$$ 
\tilde{F}= \tilde{q} \#\left( 1 - \frac{\tilde{h}}{c_0} \right) ,  
$$
$\tilde{q}$ denoting a parametrix of $\tilde{p}_2$ on the range
of $ 1 - (h/c_0)$, i.e. a symbol such that 
$$
\tilde{q} \#  \tilde{p}_2 = \tilde{p}_2 \# \tilde{q} = 1. 
$$
Thus we may think that the symbol $ \tilde{F} $ belongs to $
S^{-2}(\R^{\nu}_{t}) $.

By $ (H3) $ (see also (\ref{eq:2.2})), $ \tilde{p}_{1} $ has a parametrix.
Thus the operator $R^+_* \tilde{p}_1 R^-$ has a parametrix,
$\tilde{Q}$, which is a pseudodifferential operator in the $ y
$-variable. 
\begin{equation}\label{eq:q}
 \tilde{Q} \, \# \, R^+_* \tilde{p}_1 R^- = R^+_* \tilde{p}_1 R^- \,
 \# \, \tilde{Q}=  1 \, . 
\end{equation}
We observe that, $R^+_* \tilde{p}_1 R^- $ is a linear function with
respect to $ (y, \eta) $ and
$$
R^+_* \tilde{p}_1 R^- (y, \eta) = c_0 \tilde{p}_1 (y,\eta ) .
$$
Set
$$
\tilde{E} = \tilde{F}+\lambda^{1/2} \tilde{Q} \tilde{h}  -
\frac{1}{c_{0}} \tilde{h} \#  \tilde{F},
$$
$ \tilde{E} $ is an approximate parametrix of the operator
$\tilde{P}$, modulo $\lambda^{-1/2}$, i.e. we have
\begin{equation}
\label{eq:param}
\tilde{E}\# \tilde{P} = 1+ {\lambda^{-1/2}}
 \left [ \tilde{F}  - \frac{1}{c_{0}} \tilde{h} \#  \tilde{F}
   \right ] \tilde{p}_1.
\end{equation}
The symbols appearing in the above formula are real analytic symbols
in the classes $ S^{m}(\R^{\nu}_{t}; S^{\ell}(\R_{y})) $. 

It might be worth saying that the above identity for our approximate
parametrix relies on a couple of identities:
$$ 
\tilde{h} \# \tilde{p}_{2} = 0
$$
and 
$$ 
\frac{1}{c_{0}} \tilde{h} \# \tilde{h} = \tilde{h}.
$$
The first is an easy consequence of the definition of the symbol $ \tilde{h}
$. 

The error term obtained in Formula (\ref{eq:param}) above is a symbol
in the class $ S^{-2}(\R^{\nu}_{t}; S^{1}(\R_{y})) $.

\subsection{The metaplectic FBI transform} 
\renewcommand{\theequation}{\thesubsection.\arabic{equation}}
\setcounter{equation}{0}
\setcounter{theorem}{0}
\setcounter{proposition}{0}  
\setcounter{lemma}{0}
\setcounter{corollary}{0} 
\setcounter{definition}{0}
\setcounter{remark}{0}

In the present context we use the following definition of FBI transform:
\begin{equation}
\label{eq:mfbi} 
Tu(x)=\int_{\R^{\nu+1}} e^{- \frac 12 (x-x')^2} u(x')\, dx',
\end{equation}
where $x'=(t',y')\in \R^{\nu +1}$, $x=(t,y)\in \C^{\nu +1}$ and $u$ is
e.g. a tempered distribution. 
Defining 
$$ 
\Phi(x, x') = \frac{i}{2}(x-x')^{2} \qquad \text{and}\qquad 
\phi_{0}(x) = - \sup_{x'} \im \Phi(x, x'),
$$
we find  
$$
\phi_0 (x) = \frac {| \im x |^2} 2 \, .
$$
In the sequel, we will use the notation 
$$
\phi_0 (x) = \phi_{0,1} (t) + \phi_{0,2} (y)
$$
where $\phi_{0,1} (t)=| \im t |^2/2$ and $\phi_{0,2} (y) =( \im y
)^2/2$. 
One can show that $T$ maps 
$L^2(\R^{\nu +1})$ into  $H_{\phi_0}$ (the space of entire function on
$\C^{\nu +1}$, square integrable with respect to the measure $
e^{-2\phi_0(x)} L(dx)$, where $L(dx)=\left ( \frac i2\right)^{\nu +1}
dx\wedge d\bar{x}$ is 
the Lebesgue measure in $\R^{2 (\nu +1)}$). In the sequel we will also
use the partial FBI transform with respect to the $t$-variables only
and we still denote it by $T$. It will be clear from the context
whether we are considering a partial or a global transformation.    

We have
$$ 
T \tilde{P} u = P Tu,
$$
where, as symbols,
$$ 
P \circ \mathcal{H}_{T}(x, \xi) = \tilde{P}(x, \xi)   ,
$$
and
$$ 
\mathcal{H}_{T}(x, \xi) = (x - i \xi, \xi).
$$
Actually $ T $ is associated
with the complex canonical transformation 
\begin{equation}
\label{tc}
\C^{2(\nu+1)} \ni\left(x', -  \partial_{x'} \Phi\right)
\longmapsto (x, \partial_{x} \Phi) \in \C^{2(\nu+1)}. 
\end{equation}
Henceforth we write $ \frac{\partial f}{\partial x'} $ for a real derivative,
whereas $ \partial_{x} $ denotes the complex derivative 
$\partial_{x} = (1/2) ( \partial_{\re x} - i \partial_{\im x})$.

We observe that the range of $ \mathcal{H}_{T} $ coincides with the
range of the restriction of (\ref{tc}) to $ \R^{2(\nu +1)} $ and is an
$ I $-Lagrangian manifold (that is a Lagrangian manifold for the
non-degenerate skew-symmetric form $ \im \sigma $) in $ \C^{2(\nu+1)}
$ of the type 
$$ 
\Lambda_{\phi_{0}} = \left\{ \left( x, \frac{2}{i}
    \partial_{x}\phi_{0}(x)\right) \right\} = \left\{ \left( x, - \im
    x \right) \right\}.  
$$

\subsection{The parametrix on the FBI side}
\renewcommand{\theequation}{\thesubsection.\arabic{equation}}
\setcounter{equation}{0}
\setcounter{theorem}{0}
\setcounter{proposition}{0}  
\setcounter{lemma}{0}
\setcounter{corollary}{0} 
\setcounter{definition}{0}

We have seen in (\ref{eq:param}) that, there exists a suitable
operator $\tilde{E}$ giving a parametrix of $ \tilde{P} $:
$$
\tilde{E} \# \tilde{P} =  1+ \lambda^{-1/2}
 \left [ \tilde{F}  -  \frac{1}{c_{0}} \tilde{h} \# \tilde{F}
   \right ] \tilde{p}_1.
$$

Define
$$
\tilde{F} = F \circ \mathcal{H}_T \qquad
\tilde{E} = E \circ \mathcal{H}_T \qquad
\tilde{h} = h \circ \mathcal{H}_{T}.
$$
Since the Weyl composition and the linear canonical transformation
commute, we have
\begin{equation}
\label{eq:ep}
E \# P = 1 + \lambda^{-1/2} \left[ F - \frac{1}{c_{0}} h \# F\right] p_{1}. 
\end{equation}
For an analytic symbol $ q(t, y, \tau, \eta) $ we define the
corresponding pseudodifferential operator on holomorphic functions $ u
$ as 
\begin{multline}\label{eq:QT}
q_{S,\chi}(t, y, D_{t}, D_{y}) u(t,y)=\frac 1 {(2\pi)^{\nu+1}} \iint
e^{i \langle t-t', \tau \rangle + i  (y-y') \eta }  \times 
\\
q
\left (\frac{t+t'}{2}, \frac {y+y'}2, \tau, \eta\right ) 
\chi \left (\frac{t-t'}{S}, \frac {y-y'}S\right)
u(t', y') dt'\wedge dy'\wedge d\tau\wedge d\eta.
\end{multline}
Here the integral is computed along the path 
$$
\Gamma \ \colon \xi=(\tau, \eta) \ \left\{
\begin{array}{c}
\displaystyle{\tau = \frac{2}{i}
  \partial_{t}\phi_{0,1}\left(\frac{t+t'}{2}\right)  + 
\frac i S \overline{ (t - t') } }
\\[25pt]
\displaystyle{\eta =
  \frac{2}{i}\partial_{y}\phi_{0,2}\left(\frac{y+y'}{2}\right)  + 
\frac i S \overline{ (y - y') } }
\end{array}
\right .
$$
and $\chi (t,s)=\chi_1 (t)\chi_2 (s)$, where $\chi_{1}(t)$, $
\chi_{2}(y) $ are cut-off functions equal to $1$ near the origin .

We point out that in defining the above realizations we use the fact
that the symbol $q$ can be holomorphically continued  to a
neighbourood of $\Lambda_{\phi_0}$ in $\C^{2(\nu +1)}$.

By (\ref{eq:ep}), we have  $(E \# P )_{S,\chi }=1_{S,\chi }+
\lambda^{-1/2}A_{S,\chi }$, 
where  
\begin{equation}\label{eq:A}
A=\left[ F - \frac{1}{c_{0}} h \# F\right]
p_{1}
\end{equation}
is the FBI transform of the error term belonging to $
S^{-2}(\R^{\nu}_{t}; S^{1}(\R_{y})) $ obtained in (\ref{eq:param}).

Hence 
\begin{multline}
\label{eq:parcutoff}
E_{S,\chi } P =1+[1_{S,\chi }-1]+ [E_{S,\chi } P -(E \# P )_{S,\chi }]
+ \lambda^{-1/2}A_{S,\chi }.
\end{multline}
In the next section we proceed to estimate the errors $1_{S,\chi }-1$, 
$\lambda^{-1/2} A_{S,\chi}$ and $E_{S,\chi } P -(E \# P )_{S,\chi }$.

\subsection{Estimate of the errors} 
\renewcommand{\theequation}{\thesubsection.\arabic{equation}}
\setcounter{equation}{0}
\setcounter{theorem}{0}
\setcounter{proposition}{0}  
\setcounter{lemma}{0}
\setcounter{corollary}{0} 
\setcounter{definition}{0}
\setcounter{remark}{0}

We denote by  $ \Sigma_{1}^{\C} = \{ t = \tau = 0\} $, $ \Sigma_{2}^{\C} = \{
(0, 0)\} $ the complexifications of $ \Sigma_{1} $, $
\Sigma_{2} $, the characteristic manifold of $ P $ and the
characteristic set of the first eigenvalue of $ P $ respectively).
Then, $ \Sigma_{1} = \Sigma_{1}^{\C} \cap \Lambda_{\phi_{0}} $, $
\Sigma_{2} = \Sigma_{2}^{\C} \cap \Lambda_{\phi_{0}} = \Sigma_{2}^{\C}
$.

We keep understanding that $ x = (t, y) \in \C^{\nu+1} $ and define
$$ 
d_{1}(x) = \text{distance of}\ \left( x, \frac{2}{i}
    \partial_{x}\phi_{0}(x)\right) \text{to} \ \Sigma_{1},
$$
and
$$ 
d_{2}(x) = \text{distance of} \ \left( x, \frac{2}{i}
    \partial_{x}\phi_{0}(x)\right) \text{to} \ \Sigma_{2}.
$$

Let $ \Omega \subset \C^{ \nu +1} $ be an open subset and $ u $ a
square integrable function on $ \Omega $. We need the following
norms for $ u $: 

$$ 
\| u \|^{2}_{\phi_{0}, \Omega} = \int_{\Omega} e^{-2\phi_{0}(x)}
|u(x)|^{2} L(dx),
$$
$$
\tb{u}^{2}_{\phi_{0}, \Omega} = \int_{\Omega} e^{-2\phi_{0}(x)}
\left(d_{1}^{2}(x) + 1 \right)^{2}
|u(x)|^{2} L(dx). 
$$

Let $ S \geq 1 $ be a large parameter and define 
$$
B(0, S) =  \{ (t,y) \in \C^{\nu +1}| |(t,y)| < S\}.
$$
In general we use the following notation: if $ a $ and $ b $ are
quantities depending both on $ \lambda $ and $ 
S $, we write $a\lesssim b$ for $a\leq C b$ for a suitable $C>0$
independent of $S$ and $\lambda$.

In the proofs of the present section we will use several times and
without any further mention the following elementary remark. 
\begin{remark}
Let $A$ be an operator defined by an integral kernel $K$, i.e. 
$$
Au (x) = \int_{\Gamma} K(x,x') u(x') L(dx') , 
$$
on a possibly complex domain $ \Gamma $, and set 
$$
d(x)=d_{1}^{2}(x) + 1.
$$
Then 
\begin{multline*}
\tb{Au}_{\phi_0, \Omega}^2=\|d(\cdot) Au\|_{\phi_0, \Omega}^2 
\\
\leq \iint_{\Omega\times\Gamma} |
  d(x) e^{-\phi_0 (x)+\phi_0 (x')} K(x,x') |^2 L(dx) 
  L(dx') \times 
\\
\int_{\Omega} e^{-2\phi_0 (x')} |u(x')|^2 L(dx')
\end{multline*} 
and in order to control the norm $\tb{Au}_{\phi_0, \Omega}^2$ it is enough
to estimate the norm of the reduced kernel $d(x) e^{-\phi_0
  (x)+\phi_0 (x')} K(x,x')$. 
\end{remark} 
For the proof of the next result we refer to Sj\"ostrand
\cite{Sj-83}, Equation (1.11) and the following discussion, and to
\cite{Sj-Ast}, Equation (12.45). 
\begin{lemma}
\label{l:1}
Let $u$ be a holomorphic function in $B(0,2S)$. Then, 
there exists a positive constant $C$, independent of $ u $, $ \lambda
$ and $ S $, such that 
\begin{equation}\label{eq:sj}
\tb{ 1_{S,\chi } u -u }_{\phi_0, B(0,S) }\lesssim e^{-S/C} 
\| u \|_{\phi_0, B(0,2S) }. 
\end{equation}
\end{lemma}
The above Lemma obviously takes care of the first error term in
(\ref{eq:parcutoff}). Next we estimate the third term.
\begin{lemma}
\label{l:2}
Let $u$ be a holomorphic function in $B(0,2S)$. Then,
\begin{equation}
\tb{ A_{S,\chi }u }_{\phi_0, B(0,S) }\lesssim  S^{\nu + 4} 
\| u \|_{\phi_0, B(0,2S) }.
\end{equation}
\end{lemma} 
{\sc Proof:} We provide a rough estimate of the reduced kernel 
\begin{multline*}
\left | \left (d_1^2(x)+ 1\right )
e^{-\phi_0(x)+\phi_0 (x')} e^{i\langle x-x',\xi \rangle } \right . 
\\
\left . a\left (
\frac {x+x'}2 \, , \, \xi \right ) \chi \left ( \frac {x-x'}S\right ) \det
\left ( \frac {\partial \xi }{\partial x'} \right ) \right |.
\end{multline*}
Since both $ |x| $ and $ |x'| $ are bounded by $ C S $, due to the
cutoff, we may estimate the above quantity by
\begin{multline*}
S^{2} e^{-|x-x'|^2/S} \chi \left ( \frac {x-x'}S\right )
\\
\left | a \left ( \frac {x+x'}2 \, , \, \frac 2i \partial_x \phi_0 \left (
\frac {x+x'}2 \right ) +\frac iS \overline{(x-x')} \right )\right |.
\end{multline*} 
Now, $ a $ contains terms that can be estimated by $ (1 +
d_{1}(x))^{-2} (1 + d_{2}(x)) $, like $ Fp_{1} $, or by $ (1+
d_{1}(x))^{-N} (1 + d_{2}(x)) $, like the second term in
(\ref{eq:A}). We thus obtain the bound for the reduced kernel
$$
S^3 e^{-|x-x'|^2/S} \chi \left ( \frac {x-x'}S\right )
$$
and the conclusion follows. 
\Qed 
The following lemma takes care of the  second error term in 
(\ref{eq:parcutoff}) and is due to Sj\"ostrand \cite{Sj-83} and we
sketch its proof only to make the present paper self-contained.
\begin{lemma}
\label{l:2.5}
Let $ q(x, \xi) $ be an analytic symbol and denote by $ q_{S,\chi} $
its realization defined as in (\ref{eq:QT}). Let $ P(x, D_{x}) $ be a
differential operator such that
$$ 
P(x, D_{x}) = \sum_{k \leq m} a_{j_{1}, \cdots , j_{k}} U_{j_{1}}
\circ \cdots \circ U_{j_{k}}, \qquad a_{j_{1}, \cdots , j_{k}} \in \C,
$$
where
$$ 
U_{j}(x, D_{x}) = \langle \alpha_{j}, x \rangle + \langle \beta_{j},
D_{x}\rangle,
$$
$ j = 1, \ldots, d $, $ x \in \C^{n} $. Then
$$ 
q_{S, \chi}(x, D_{x}) \circ P(x, D_{x}) = (q \# P)_{S,\chi}(x, D_{x})  +
\frac{1}{S} \sum_{j=1}^{n} R^{(j)}_{S,\partial_{x_{j}}\chi}(x, D_{x}). 
$$
Here $ R^{(j)} $ are analytic symbols whose realization is defined as in
(\ref{eq:QT}), replacing the cutoff $ \chi $ with $ \partial_{x_{j}}\chi $.
\end{lemma}
{\sc Proof:}
The proof boils down to computing $ q_{S,\chi}(x, D_{x}) \circ
D_{x_{k}}  $ and $ q_{S, \chi}(x, D_{x}) \circ x_{k}  $. For a
holomorphic function $ v $ we have
\begin{multline*}
\frac 1i d_{(x', \xi )} \left (  q\left ( \frac {x+x'}2 , \xi
  \right) e^{i\langle x-x', \xi \rangle }  \chi \left ( \frac {x-x'} S
  \right )  v(x') \right. 
\\ 
\left. \vphantom{\frac {x-x} S}
(-1)^{k-1} dx'_{1}\wedge \cdots \wedge \widehat{dx'}_{k}
  \wedge \cdots \wedge dx_{n}' \wedge d\xi  \right)
\\
= D_{x'_{k}} \left( q\left ( \frac {x+x'} 2 , \xi\right
  )  e^{i\langle x-x', \xi \rangle } \right) \chi\left(\frac{x-x'}S
  \right)  v(x')\,  dx'\wedge d\xi
\\
+  q\left ( \frac {x+x'} 2 , \xi \right
  ) e^{i\langle x-x', \xi \rangle }  \chi \left ( \frac {x-x'} S
  \right )  D_{x'_{k}} v(x')\, dx'\wedge d\xi
\\
+ q\left ( \frac {x+x'}2 , \xi\right)
  e^{i\langle x-x', \xi \rangle } v(x' ) 
\\
\frac 1S
\left( D_{x'_{k}} \chi \, dx' \wedge d\xi +  D_{\bar{x}'_{k}} \chi
  dx'_{1}\wedge \cdots \wedge d\bar{x}'_{k} 
  \wedge \cdots \wedge dx_{n}' \wedge d\xi \right).
\end{multline*}
Using Stokes formula we find
\begin{multline*}
q_{S, \chi } D_{x_{k}} v 
\\
= -\int_{\Gamma} D_{x'_{k}} \left( q\left (
    \frac {x+x'}2 , \xi \right)
  e^{i\langle x-x', \xi \rangle } \right) \chi\left ( \frac {x-x'} S
  \right )  v(x')  dx'\wedge d\xi
\\[3pt]
- \frac 1S \int_{\Gamma_x} q\left ( \frac {x+x'}2 , \xi\right
  ) e^{i\langle x-x', \xi \rangle } v(x' )  
\\
\left( D_{x'_{k}} \chi \, dx' \wedge d\xi +  D_{\bar{x}'_{k}} \chi
  dx'_{1}\wedge \cdots \wedge d\bar{x}'_{k} 
  \wedge \cdots \wedge dx_{n}' \wedge d\xi \right)
\\
= (q\# \xi_{k} )_{S, \chi } v + \frac 1S q_{S,\partial_{x_{k}}\chi} v
+ \frac{1}{S} q_{S, \bar{\partial}_{x_{k}} \chi} v.
\end{multline*}
Here the third term above has a volume form slightly different from
that in (\ref{eq:QT}).

On the other hand we have
\begin{multline*}
q_{S, \chi } x_k v
\\
=\int_{\Gamma}
q\left ( \frac {x+x'} 2 , \xi  \right) 
\left( \frac {x_k+x_k'}{2}\right ) e^{i\langle x-x', \xi \rangle }
\chi \left ( \frac {x-x'} S
  \right )  v(x')\,  dx'\wedge d\xi
\\
-\int_{\Gamma}
q\left ( \frac {x+x'} 2 , \xi  \right)  
\left ( \frac {x_k-x_k'}{2}\right ) e^{i\langle x-x', \xi \rangle }
\chi \left( \frac {x-x'} S\right )  v(x')\,  dx'\wedge d\xi .
\end{multline*}
Noticing that
\begin{multline*}
\int_{\Gamma}
q\left( \frac {x+x'} 2 , \xi  \right)  \left(\frac
{x_k-x_k'}{2}\right) e^{i\langle x-x', \xi \rangle }
\chi\left( \frac {x-x'}S\right )  v(x')\,  dx'\wedge d\xi
\\
=\int_{\Gamma}
q\left( \frac {x+x'} 2 , \xi  \right)  \left(\frac {D_{\xi_k}
}{2}\right) e^{i\langle x-x', \xi \rangle} 
\chi\left( \frac {x-x'}S\right) v(x')\,  dx'\wedge d\xi
\end{multline*}
and arguing as above we deduce that 
$$
q_{S, \chi } x_k v=(q\# x_k)_{S,\chi} v+ \frac {1}{2S} \sum_{j=1}^n 
q_{S, \bar{\partial}_{x_{j}} \chi} v.
$$
We point out that the integral defining the last term above is taken
with respect to the measure $(-1)^{n+k-1} d\bar{x}_j\wedge dx'\wedge
d\xi_1'\wedge\cdots \wedge \widehat{d\xi'}_k\wedge \cdots \wedge
d\xi_n'$.   
The conclusion follows by iteration. 
\Qed

\begin{lemma}
\label{l:3} 
Let $u$ be a holomorphic function in $B(0,2S)$. Then, there exists a
positive constant $C$ independent of $u$, $\lambda$ and $S$, such that  
\begin{equation}\label{eq:let}
\tb{E_{S,\chi }P u -(E \# P )_{S,\chi }u }_{\phi_0, B(0,S) }\lesssim 
e^{-S/C}  \| u \|_{\phi_0, B(0,2S) }. 
\end{equation} 
\end{lemma} 
{\sc Proof:} By Lemma~\ref{l:2.5} the proof reduces to estimate a term
of the form $\tb{ R^{(j)}_{S,\partial_{x_{j}}\chi}u }_{\phi_0, B(0,S) }$.
Since the support of the function $\partial_{x_{j}}\chi$ is away from
the origin we get the following bound for the reduced kernel 
\begin{multline*}
\left | \left (d_1^2(x)+ 1 \right )
e^{-\phi_0(x)+\phi_0 (x')} e^{i\langle x-x',\xi \rangle } \right . 
\\
\left . a\left (
\frac {x+x'}2 \, , \, \xi \right )\partial_{x_{j}} 
\chi\left(\frac{x-x'}S\right ) 
\det \left(\frac {\partial \xi }{\partial x'}\right) \right |
\\
\lesssim  S^2
e^{-S/C} \partial_{x_{j}} \chi \left ( \frac {x-x'}S\right )
\end{multline*}
and the conclusion follows. 
\Qed

The next estimate is a direct consequence of Lemma~\ref{l:1}, Lemma
\ref{l:2} and Lemma~\ref{l:3}. 
\begin{proposition}
For every  $u$ holomorphic in $B(0,2S)$ we have
\begin{equation}\label{eq:ester}
\tb{E_{S,\chi } P u - u }_{\phi_0, B(0,S) }\lesssim 
\gamma (S,\lambda )\| u \|_{\phi_0, B(0,2S) }
\end{equation}
where
\begin{equation}
\label{eq:gamma}
\gamma (S,\lambda ) =  e^{-S/C}+\lambda^{-1/2} S^{\nu + 4}. 
\end{equation}
\end{proposition} 

\subsection{The a priori estimate} 
\renewcommand{\theequation}{\thesubsection.\arabic{equation}}
\setcounter{equation}{0}
\setcounter{theorem}{0}
\setcounter{proposition}{0}  
\setcounter{lemma}{0}
\setcounter{corollary}{0} 
\setcounter{definition}{0}

In this section we prove an a priori estimate for the localized
operator. We start by estimating the the action of the parametrix $
E_{S, \chi} $  between the spaces defined previously.

\begin{lemma}
\label{l:4} 
We have the estimate
\begin{equation}
\label{eq:estE}
\tb{ E_{S,\chi } u}_{\phi_0 , B(0,S)} \lesssim
\lambda^{1/2} \| u\|_{\phi_0 , B(0,2S)}. 
\end{equation}
for every $u$ holomorphic in $B(0,2S)$. 
\end{lemma}
{\sc Proof:} 
Recalling that 
$$
E=F+\lambda^{1/2} Q h - \frac{1}{c_{0}} h\# F.
$$
the proof of (\ref{eq:estE}) reduces to estimate the above three terms. 
First we show that 
\begin{equation}\label{eq:estF}
\tb{F_{S, \chi } u }_{\phi_0 , B(0,S)}\lesssim  \|  u \|_{\phi_0 , B(0,2S)}\, .
\end{equation} 
The idea of the proof is due to Sj\"ostrand, \cite{Sj-83}.

We split the integral in the left hand side of (\ref{eq:estF}) into
two region: 
$$
(1)\quad x'\, :\quad  | x-x'|\le \frac SC \qquad \text{and}\qquad (2) \quad
x'\, :\quad | x-x'|>
\frac SC \, .
$$ 
Hence, we write 
$$
\tb{F_{S, \chi } u }_{\phi_0 , B(0,S)}=
\tb{F^{(1)}_{S, \chi } u}_{\phi_0 , B(0,S)}+
\tb{F^{(2)}_{S, \chi } u }_{\phi_0 , B(0,S)}
$$
where the symbol $F^{(i)}_{S, \chi }$ means that the
integral in the norm is taken in the region $(i)$.

In $(2)$ we have that the reduced kernel can be estimated
by $ S^{2} e^{-S/C}$. Hence, possibly taking a smaller $C$ (independent
of $u$, $S$ and $\lambda$), we get 
$$
\tb{ F^{(2)}_{S,\chi } u}_{\phi_0 , B(0,S)}\lesssim  \| u\|_{\phi_0 , B(0,2S)}\, .
$$
Now, we observe that, taking the constant $C$ large enough, we have 
$ \chi  \left (  \frac {x-x'}
  S\right )=1$  in the region $(1)$.  

In order to estimate $\tb{ F^{(1)}_{S,\chi } u}_{\phi_0 ,
    B(0,S)}$ first we replace the contour $\Gamma \cap \{ |
x-x'|\le S/C\}$, defined in (\ref{eq:QT}), with the singular contour
$\tilde{\Gamma}\cap \{ |x-x'|\le S/C\}$ (see \cite{Sj-Ast}), where 
\begin{equation}
\label{eq:gammat}
\tilde{\Gamma}\, :\quad \xi =(\tau, \eta) = \frac 2i \partial_x \phi_0
\left (\frac {x+x'}2
\right ) + \frac iC \frac {\overline{(x-x')}} {|x-x'|}\, .
\end{equation}
Here we used again the notation $ x=(t, y) $, $ x'=(t', y') $.
The reduced kernel of $F^{(1)}_{S,\chi }$ can be estimated
(modulo constants) by 
\begin{multline*}
K(x,x')=
\left(d_{1}^{2}(x) + 1 \right) \left |
  F \left( \frac {t+t'} 2 , \tau \right ) \right | e^{-|x-x'|/C} 
(1+ |x-x'|^{-\nu -1} ),
\end{multline*}
where $\tau $ is given by (\ref{eq:gammat}). 
We point out that $ \tau $ as defined in (\ref{eq:gammat}) is a
function of $ y $ and $ y' $, however the following estimate holds
\begin{multline*}
\left | F \left ( \frac {t+t'} 2 , \tau \right ) \right | \leq C_{1} \left
  (1+d_1\left ( \frac {t+t'}2 \right ) \right )^{-2} \\
\leq C_{2} (1 + | t-t'|)^2 ( 1 +d_1 (t))^{-2}, 
\end{multline*}
where $ C_{1} $ and $ C_{2} $ are positive constants depending on $ F
$ and $ C $ only.
We have that 
$$
K(x,x')\lesssim (1 + | x-x'|)^2 e^{-|x-x'|/C} (1+ |x-x'|^{-\nu -1} ) 
$$
and we conclude 
$$
\tb{F^{(1)}_{S, \chi } u }_{\phi_0 , B(0,S)}\lesssim \|  u \|_{\phi_0 , B(0,2S)}
$$
hence (\ref{eq:estF}) follows. 
Repeating the same kind of argument as above and using the fact that
for every $N$ there exists a positive constant $C_N$ such that 
$$
|(Q  h)(x,\xi )|\lesssim C_N \frac 1{1+d_2(x)} (1+d_1(x))^{-N} 
$$
we deduce that 
\begin{equation}\label{eq:estQh}
\tb{(Q h)_{S, \chi } u }_{\phi_0 , B(0,S)}\lesssim  
\| u \|_{\phi_0 , B(0,2S)}\, .
\end{equation} 
Moreover, the above arguments and the estimate 
$$
\left| \left(\frac{1}{c_{0}}  h\# F\right) (x,\xi )\right|\lesssim C_N
(1+d_1(x))^{-N}
$$ 
yield that 
\begin{equation}\label{eq:qhpf}
\tb{ c_{0}^{-1} ( h\# F)_{S,\chi } u}_{\phi_0 , B(0,S)}\lesssim  
\| u \|_{\phi_0 , B(0,2S)}\, .
\end{equation} 
Then (\ref{eq:estE}) follows by (\ref{eq:estF}),
(\ref{eq:estQh}) and (\ref{eq:qhpf}).
\Qed 
\begin{remark}
\label{rem:perdita}
We point out that the factor $ \lambda^{1/2} $ in front of the norm in
(\ref{eq:estE}) is related to the fact that the operator $ P $ has a
loss of $ 3/2 $ derivatives (see \cite{Helffer1})
\end{remark}

\begin{proposition}
\label{p:2}
We have the estimate
$$
\tb{u}_{\phi_{0}, B(0, S)}  
\lesssim \lambda^{1/2} \| P u \|_{\phi_{0},B(0, 2S)} 
+ \gamma(S, \lambda) \|u\|_{\phi_{0}, B(0, 2S)},
$$
for every holomorphic function $ u $ on $ B(0, 2S) $. 
Here $ \gamma(S, \lambda) $ is given by (\ref{eq:gamma}).
\end{proposition}
{\sc Proof:}
Lemma~\ref{l:4} yields that 
$$
\tb{E_{S,\chi }v}_{\phi_{0}, B(0, S)}  
\lesssim \lambda^{1/2} \| v \|_{\phi_{0},B(0, 2S)}
$$
for every $v$ holomorphic in $B(0,2S)$. 
Hence, taking $v=Pu$, 
\begin{multline*}
 \sqrt{\lambda}\| Pu  \|_{\phi_{0}, B(0, 2S)}\gtrsim 
\tb{E_{S,\chi }Pu}_{\phi_{0},B(0, S)}\geq 
\\
\tb{u}_{\phi_{0},B(0, S)}  
-\tb{ E_{S,\chi }Pu-u}_{\phi_{0},B(0, S)} 
\end{multline*}
and the conclusion follows by (\ref{eq:sj}) and (\ref{eq:ester}).\Qed

\section{The local a priori estimate}
\renewcommand{\theequation}{\thesection.\arabic{equation}}
\setcounter{equation}{0}
\setcounter{theorem}{0}
\setcounter{proposition}{0}  
\setcounter{lemma}{0}
\setcounter{corollary}{0} 
\setcounter{definition}{0}
\setcounter{remark}{0}
%
\renewcommand{\thetheorem}{\thesection.\arabic{theorem}}
\renewcommand{\theproposition}{\thesection.\arabic{proposition}}
\renewcommand{\thelemma}{\thesection.\arabic{lemma}}
\renewcommand{\thedefinition}{\thesection.\arabic{definition}}
\renewcommand{\thecorollary}{\thesection.\arabic{corollary}}
\renewcommand{\theequation}{\thesection.\arabic{equation}}
\renewcommand{\theremark}{\thesection.\arabic{remark}}

The purpose of this section is to provide local \textit{a priori}
estimates which will allow us to deduce a theorem on the propagation
of the regularity. From now on our ambient space is $ n
$-dimensional.

Let $\lambda \geq 1$ be a large parameter and denote by 
$$ 
D = \frac{1}{i} \partial\, , \qquad \tilde{D} = \frac{1}{\lambda} D\, .
$$
We use the $ \lambda $-Fourier transform:
$$ 
\hat{u}(\xi) = \int e^{-i\lambda x \xi} u(x) dx,
$$
$$
 u(x) =
\left(\frac{\lambda}{2\pi}\right)^{n} \int e^{i\lambda x \xi}
\hat{u}(\xi) d\xi.
$$
In this setting we use the following definition of FBI transform:
$$ 
Tu(x, \lambda) = \int_{\R^{n}} e^{- \frac{\lambda}{2}(x - x')^{2}}
u(x') dx',
$$
$ \lambda \geq 1 $. $T$ is associated with the complex canonical
transformation $\mathcal{H}_T \, (x,\xi )\mapsto (x-i\xi , \xi )$. 
Moreover, we have  
$$
\mathcal{H}_T (\R^{2n}) = \Lambda_{\phi_0}=\{
(x,(2/i)\partial_x \phi_0(x))| x\in\C^n\},
$$ 
with $\phi_0(x)=(1/2) |\im x|^2$.   

If $ u \in \mathcal{S}' (\R^n )$ then $ Tu $ is a holomorphic function
of  $ x \in \C^{n}$ and moreover if $ u $ belongs to $ L^{2}(\R^{n}) $
then $ Tu \in L^{2}(\C^{n}, e^{-2\lambda \phi_{0} (x)} L(dx))$.

We recall the characterization of the analytic wave front set in the
FBI setting (see e.g. \cite{Sj-Ast}): a point $(x_{0}, \xi_{0}) \in
\R^{2n}$  does not belong to $ WF_{a}(u)$ iff there exist a positive
$ \epsilon $, a neighborhood $ V $ of $ 
x_{0} - i \xi_{0} $ in $ \C^{n} $ and a positive constant $ C_{V} $
depending on $ V $, such that
$$ 
| e^{-\lambda \phi_{0}(x)} Tu(x, \lambda) | \leq C_{V} e^{-\epsilon \lambda}
$$
uniformly in $ V $ for $ \lambda $ large enough.
When on the FBI side, we write $ x_{0} $ instead of $ x_{0} - i \xi_{0}$.

Denoting again by 
$ \tilde{P} $ the given operator before the FBI transform, we have
\begin{equation}
\label{eq:pt} 
\lambda^{-2} \tilde{P}(x, D) = \tilde{p}_{2}(x, \tilde{D}) +
\lambda^{-1} \tilde{p}_{1}(x, \tilde{D}) + \mathcal{O}(\lambda^{-2}),
\end{equation}
for $x$ in an open set of $\R^n$. 

Let $\tilde{\Sigma}_1$,  
$\tilde{\Sigma}_2$ the real analytic manifolds of Assumptions
$(H1)$ and $(H3)$ of Section~\ref{s:sr}, respectively.
It is always possible to perform a homogeneous canonical
transformation such that 
$$ 
\tilde{\Sigma}_{1} = \{ t = \tau = 0 \},
$$
$$ 
\tilde{\Sigma}_{2} = \{ t = \tau = 0, y = \eta = 0 \},
$$
where, after the canonical transformation the new variables have been
written as $ x = (t, y, s) \in \R^{n_{t}}\times \R \times \R^{n_{s}}$,
where $ n = n_{t} + 1 + n_{s} $.

Let $ (z, \zeta) = (0, \bar{y}, \bar{s}; 0, \bar{\eta}, \bar{\sigma})
$ be a point in $ \tilde{\Sigma}_{1} \setminus \tilde{\Sigma}_{2} $,
then we denote by $(\tilde{P})_{(z,\zeta)}$ the localization at $ (z,
\zeta) $ of the symbol $ \tilde{p}(x, \xi) $.
We have that 
\begin{eqnarray*}
(\tilde{P})_{(z,\zeta)}(x,\xi , \lambda ) &=& \sum_{k=0}^{1} \sum_{|
  \alpha +\beta |=2 - 2k} \frac
{\lambda^{-k}}{\alpha ! \beta !} \tilde{p}_{2-k (\beta )}^{(\alpha )} (z,\zeta )
(x-z)^\beta (\xi -\zeta )^\alpha \\
 & = &  \sum_{k=0}^{1} \sum_{|
  \alpha +\beta |=2 - 2k} \frac
{\lambda^{-k}}{\alpha ! \beta !} \tilde{p}_{2-k (\beta )}^{(\alpha )}
(0, \bar{y}, \bar{s}; 0, \bar{\eta}, \bar{\sigma}) t^{\beta} \tau^{\alpha}.
\end{eqnarray*}
Using the natural homogeneity we denote also by
$[\tilde{P}]_{(z,\zeta)}$ the symbol
\begin{eqnarray}
\label{eq:loc22}
[\tilde{P}]_{(z,\zeta)} & = & \sum_{k=0}^{1} \sum_{|
  \alpha +\beta |=2 - 2k} \frac{1}{\alpha ! \beta !} \tilde{p}_{2-k
  (\beta )}^{(\alpha )} (z,\zeta ) (\delta x)^{\beta}
(\delta\xi)^{\alpha} \nonumber 
\\
& = & \sum_{k=0}^{1} \sum_{|
  \alpha +\beta |=2 - 2k} \frac{1}{\alpha ! \beta !} \tilde{p}_{2-k
  (\beta )}^{(\alpha )}(0, \bar{y}, \bar{s}; 0, \bar{\eta},
\bar{\sigma}) (\delta t)^{\beta} (\delta \tau)^{\alpha}. 
\end{eqnarray}
On the other hand, for $(z,\zeta )\in \tilde{\Sigma}_2$, ie $ (z,
\zeta) = (0, 0, \bar{s}; 0, 0, \bar{\sigma}) $, we define
\begin{eqnarray*}
(\tilde{P})_{(z,\zeta)}(x,\xi , \lambda ) &=& 
 \sum_{|\alpha +\beta |=2} \frac {1}{\alpha ! \beta !}
\tilde{p}_{2\, (\beta )}^{(\alpha )}(0, 0, \bar{s}; 0, 0,
\bar{\sigma}) t^{\beta} \tau^{\alpha} \\
&  &
+ \lambda^{-1} 
\left[-\mu_{0} + \tilde{\ell}_{1}(y, \eta) 
  \right].
\end{eqnarray*} 
Here $ \mu_{0} $ has been defined in (\ref{eq:2.3}), while $
\tilde{\ell}_{1} $ is a  complex linear form.

Correspondingly we also define
\begin{multline*}
[\tilde{P}]_{(z,\zeta)}(x,\xi )= \sum_{| \alpha +\beta |=2} \frac
1{\alpha ! \beta !} \tilde{p}_{2\, (\beta )}^{(\alpha )} (z,\zeta )
(\delta t)^\beta (\delta\tau)^\alpha -\mu_0
\\
+ \lambda^{-1/2} \tilde{p}_{1}(\delta y,  \delta \eta),
\end{multline*}
where $ \tilde{p}_{1} $ has been defined in
(\ref{eq:2.3bis}).

One can show that there exists a unique formal classical analytic
pseudodifferential operator of order $0$, $P$, such that 
$$ 
T \tilde{P} u = P T u.
$$
We write 
$$
P(x,\xi,\lambda)=p_2(x,\xi )+\lambda^{-1} p_1(x,\xi )+
\mathcal{O}(\lambda^{-2}).  
$$
We denote by $\Sigma_1$, $\Sigma_2$ the manifolds
$\mathcal{H}_T(\tilde{\Sigma}_j)$, $j=1,2$; $\Sigma_j \subset
\Lambda_{\phi_0}$. We also denote by $\Sigma_j^{\C}$, $j=1,2$, the
complexifications of $\Sigma_j$. For $(z,\zeta )\in \Sigma_1 \setminus
\Sigma_{2}$ or in $ \Sigma_{2} $, we define $(P)_{(z,\zeta )}$ and
$[P]_{(z,\zeta )}$ as above. One can show
that, for the localized operators the relation
$[P]_{(z,\zeta )}T=T[\tilde{P}]_{(z,\zeta )}$ holds (in this formula $T$
stands for the metaplectic FBI transform given in (\ref{eq:mfbi})). 
Since $ P $ is a pseudodifferential operator, we must discuss
its action on spaces of the type
$  L^{2}(\Omega, e^{-2\lambda\phi(x)} L(dx)) $, $\Omega\subset \C^n$,
open set, $x_0\in \Omega$, for  a suitable weight function $ \phi $. 

In the sequel it will be useful to deal with a
\textit{deformation} $ \phi $ of $ \phi_0$ (see \cite{Sj-83}). 
Let $ W $ be a complex neighborhood of $ (x_{0}, \xi_{0}) $ in $
\C^{2n} $, such that $ W \cap \Lambda_{\phi_{0}} $ is a suitably small
neighborhood of $ (x_{0}, \xi_{0}) $ in $ \Lambda_{\phi_{0}} $ and let 
\begin{equation}
\label{eq:F2}
F\colon W \rightarrow \C^{2n}
\end{equation}
be a $ C^{\omega} $ map. We assume that $ F $ satisfies the conditions 
\par\medskip
\begin{itemize}
\item[$(a)$] $ F $ is close to the identity map e.g. in the $ C^{1}(W)
  $ norm.  
\par\noindent
One can show that, since $ F $ is close to the identity map, $ F(W\cap
\Lambda_{\phi_{0}})$ has an injective projection onto $ \C_{x}^{n}
$. Thus it is a graph.
\item[$(b)$]
There exists a real valued non negative function 
$$ 
\phi \in
C^{\infty}(\pi_{x}(F(W\cap \Lambda_{\phi_{0}}))) 
$$ 
such that
$$ 
F(W\cap\Lambda_{\phi_{0}}) = \Lambda_{\phi} \cap F(W).
$$
($ \pi_{x} \colon \C^{n}\times\C^{n} \rightarrow \C^{n} $ denotes the
projection onto the first factor.)
\item[$(c)$]
For $ j=1,2 $ we have $\Lambda_{\phi}\cap \Sigma^{\C}_{j} \cap F(W) =
F(\Sigma_{j}\cap W)$.
\end{itemize}
\par\medskip
We shall actually construct $ F $ as the flow out of a suitable
Hamiltonian field tangent to $ \Sigma_{j}^{\C} $ at points of $
\Sigma_{j}^{\C} $, $ j = 1,2 $.

Since $ F $ is close to the identity map, because of Assumptions
$(H1)$, $(H3)$, we have
$$ 
{p_{2}}_{|_{\Lambda_{\phi}}} \sim \dist^{2}\left(\left(x,
    \frac{2}{i}\partial_{x}\phi(x)\right); \Sigma_{1}^{\C} \cap
  \Lambda_{\phi}\right) = d_{1\phi}^{2}(x),
$$
$$ 
\left(p^{s} + \tr^{+}p_{2}\right)_{|_{\Sigma_{1}}} 
 \sim \dist\left(\left(x,
    \frac{2}{i}\partial_{x}\phi(x)\right); \Sigma_{2}^{\C} \cap
  \Lambda_{\phi}\right) = d_{2\phi}(x).
$$

\medskip\noindent 
Let $ \Omega \subset \C^{n} $, $ \Omega_{1} \subset\!\subset \Omega $,
$ x \in \Omega_{1} $. We define
$$
Pu(x, \lambda) 
= \left(\frac{\lambda}{2\pi}\right)^{n} \iint
e^{i\lambda (x-x')\xi} P (x, \xi, \lambda) u(x') dx' \wedge d\xi,
$$
where the integration is performed along the contour
\begin{equation}
\label{eq:con} 
\xi = \frac{2}{i} \partial_{x}\phi(x) + i K (\overline{x-x'}), \qquad
|x-x'| \leq r,
\end{equation}
where $ \phi $ is a phase function satisfying the above hypotheses, $
r $ is a small positive constant such that $ \dist(\Omega_{1},
\complement\Omega) > r $ and $ K $ is a positive constant large enough
so that, when $\xi$ is in the contour defined in (\ref{eq:con}), we
have the inequality
\begin{equation}\label{eq:sal}
e^{-\lambda (\phi (x)-\phi (x'))} \left | e^{i\lambda (x-x')\xi
}\right | \leq e^{-\lambda K  |x-x'|^2/C},
\end{equation}
for some $C>0$. Furthermore $ K $ and $ r $ are such that the contour
(\ref{eq:con}) is contained in the open set $ W \ni (x_{0}, \xi_{0}) $.

\begin{remark}
\label{rem:K}
In what follows we shall need to absorb a number of error terms and
this will be done by choosing $ K $ large enough. More precisely the
size of $ K $ will depend on $ S $ as well as on a number of constants
depending only on the given operator $ P $. On the other hand it
suffices to choose  $ S \leq S_{0} $, where $ S_{0} $ is a fixed
positive quantity depending only on the data. It turns out that the
contour (\ref{eq:con}) is contained in $ W $ provided that $ r $ is
small enough depending on the problem's data.
\end{remark}

The above realization allows us to prove the continuity of $ P $ between the
function spaces $ L^{2}_{\phi}(\Omega) $ and $ L^{2,2}_{\phi}(\Omega)
$. 
Here $ L^{2}_{\phi}(\Omega) $ is the set of all locally square
integrable functions defined on $\Omega$ equipped with the norm  
$$ 
\|u\|^{2}_{\phi, \Omega} = \int_{\Omega} e^{-2\lambda \phi(x)}
|u(x)|^{2} L(dx)
$$
and $ L^{2,2}_{\phi}(\Omega) $ defined by the norm
$$
\tb{u}^{2}_{\phi, \Omega} = 
\int_{\Omega} e^{-2\lambda \phi(x)} 
(d_{1\phi}^{2}(x) + \lambda^{-1})^{2}
|u(x)|^{2} L(dx).
$$

\medskip\noindent 

Next we establish a relation between the norms used in Section
\ref{s:3} for the localized operators and the norms of the present
Section on a small ball centered at points of $ \Sigma_{2}^{\C} \cap
\Lambda_{\phi} $.

\begin{lemma}
Let $F$ be the map defined in (\ref{eq:F2}). Then $ F $ is close to
the identity map in the $C^1$-norm, i.e. $\|
F-I\|_{C^1}=\mathcal{O}(\varepsilon)$, $ \epsilon $ a small positive
parameter. Let $ z = (0, 0, s^{*}) $ be a fixed point in $
\Sigma_{2}^{\C} \cap \Lambda_{\phi} $. Denote by 
$$ 
V(z, S, \lambda) = B_{(t, y)}((0, 0), S \lambda^{-1/2}) \times
B(s^{*}, S \lambda^{-1/2}).
$$
Then, for every $u$ holomorphic in $V(z , S,\lambda )$, 
\begin{multline}
\label{eq:4.5}
\| u\|^2_{\phi , V(z , S, \lambda )}\\
= 
\lambda^{-n} (1+\mathcal{O}_S(\varepsilon )) e^{-2\lambda \phi(z)} \int_{B(0,S)}
e^{-2\phi_{0}''(s)} \|v(\cdot , \cdot , s 
)\|^2_{\phi_{0}',B_{t,y} (0,S)}  \, L(ds)      
\end{multline}
where
$$ 
\phi_{0}'(t, y) = \frac{| \im (t, y)|^{2}}{2}, \qquad \phi_{0}''(s) =
\frac{| \im s|^{2}}{2}
$$
and
\begin{equation}
\label{eq:U}
v(t,y,s) e^{i \lambda^{1/2} \langle x', \zeta\rangle } =
u\left(\frac{t}{\lambda^{1/2}}, \frac{y}{\lambda^{1/2}}, 
\frac{s-s^{*}}{\lambda^{1/2}}\right). 
\end{equation}
Here the factor $ e^{i \lambda^{1/2} \langle x', \zeta\rangle } $ is
just the function $ e^{i \lambda \langle x-z, \zeta\rangle} $ in the new
coordinates. Moreover
\begin{multline}
\label{eq:4.7}
\tb{u}_{\phi , V(z,S,\lambda )}^2
\\
=\lambda^{-(n + 2)} (1+\mathcal{O}_S(\varepsilon )) e^{- 2 \lambda \phi(z)}
\int_{B(0,S)} e^{-2 \phi_{0}''(s)}
\tb{v(\cdot , \cdot , s 
)}^2_{\phi_{0}',B_{t,y} (0,S)}  \, L(ds).
\end{multline}
 Here the symbol $\mathcal{O}_S(\varepsilon)$ denotes a quantity $
 \mathcal{O}(\epsilon) $ such that 
$\mathcal{O}(\varepsilon)/\varepsilon $ has a polynomial bound with
respect to $S$ for $ \epsilon $ small.
\end{lemma} 
{\sc Proof:} 
We define 
$$
\Phi (x')=\lambda [\phi (x)-\phi (z)-\nabla \phi (z) (x-z)], 
$$
with $x' = \lambda^{1/2} (x-z)$, $ x = (t, y, s) $ and $ z $ defined
above. 

\medskip \noindent 
{\bf Step $1$:} We want to show that  
\begin{equation}
\label{eq:st1}
\| u\|^{2}_{\phi , V(z, S, \lambda) }= \lambda^{-n} e^{-2\lambda
  \phi(z)} (1 + \mathcal{O}_S (\epsilon )) \| v\|^{2}_{\phi_0, B_{(t,
    y)}(0,S)\times B_{s}(0, S)}. 
\end{equation}
Let us write $ x = z + \lambda^{-1/2} x' $. Then
$$ 
\tilde{D}_{x} = \lambda^{-1/2} D_{x'} \qquad x-z = \lambda^{-1/2} x'
$$
where $ x-z $ and $ x' $ in the latter equation are meant to be
multiplication operators.

We denote by $ P_{(z, \zeta)} $ the operator $ (P)_{(z, \zeta)}(x,
\tilde{D}_{x}) $ realized as a differential operator.

Thus we have the following relation between differential polynomials 
\begin{equation}
\label{eq:UU}
P_{(z,\zeta )}(x, \tilde{D}_{x}) = \lambda^{-1} [P]_{(z,\zeta )}(x',
D_{x'}). 
\end{equation}
We get 
$$
|u(x)|e^{-\lambda \phi (x)} =|v(x')|
e^{-\Phi (x') -\lambda \phi (z)}.
$$ 
Hence 
$$
\| u\|_{\phi , V(z , S, \lambda) }= \lambda^{-n/2} e^{- \lambda
  \phi(z)} \| v\|_{\Phi, B_{(t, y)}(0,S)\times B_{s}(0, S)}.  
$$ 
If the Lipschitz norm of $F-\id $ is bounded by $\mathcal{O} (\epsilon
)$ then  
\begin{multline*}
\pi_{\zeta}\left( (F-\id)(x, \frac{2}{i}\partial_{x}\phi_{0}(x)) -
  (F-\id)(z, \frac{2}{i}\partial_{x}\phi_{0}(z)) \right) \\
=
\frac{2}{i} \left[\partial_{x}\phi (x)- \partial_{x} \phi_0 (x)
\right] - \frac{2}{i} \left[\partial_{x} \phi (z)- \partial_{x} \phi_0
  (z) \right]
=\mathcal{O} (\epsilon )(x-z),
\end{multline*}
i.e. $\nabla \phi (x)-\nabla \phi (z)=[\nabla^2 \phi_0+
\mathcal{O} (\epsilon )](x-z)$. It follows that 
\begin{multline*}
\Phi (x') = \lambda \left[ \int_{0}^{1}  \langle \nabla \phi(z + \rho
  (x-z)) - \nabla \phi(z), x-z\rangle d\rho \right]
\\
=\phi_0(x')+ \mathcal{O} (\epsilon ) |x'|^2. 
\end{multline*}
Then, we deduce 
$$
\| v\|_{\Phi, B_{(t, y)}(0,S)\times B_{s}(0, S)}= (1 + \mathcal{O}_S
(\epsilon )) \| v\|_{\phi_0, B_{(t, y)}(0,S)\times B_{s}(0, S)}.
$$
Step $1$ is completed.  

\medskip \noindent
{\bf Step $2$:} We want to show that 
\begin{equation}
\label{eq:ttb}
\tb{u}_{\phi, V(z, S, \lambda)}=\lambda^{-1-\frac n2} e^{-\lambda \phi(z)}
  (1 + \mathcal{O}_S (\epsilon ))\tb{v}_{\phi_0, B_{(t, y)}(0,S)\times B_{s}(0, S)}.  
\end{equation}
Since
\begin{multline*}
\partial_{t} \phi(x) = \mathcal{O}(\epsilon) (x-z)_{t} +
\partial_{t}\phi_{0}(t)
= \mathcal{O}(\epsilon) t + \partial_{t}\phi_{0}(t) \\
= \lambda^{-1/2} \left( \mathcal{O}(\epsilon) t' +
  \partial_{t}\phi_{0}(t') \right), 
\end{multline*}
we have that
$$ 
d_{1\phi}(x) \sim \frac{1}{\lambda} \left( (1+\mathcal{O}(\epsilon))
  |t'|^{2} + |\partial_{t}\phi_{0}(t')|^{2}\right)^{1/2} \sim
\frac{1}{\lambda} d_{1}(x').
$$
Then, we get that Formula (\ref{eq:ttb}) holds.  \Qed 

The next proposition is the core estimate of the present
section. Actually the microlocal regularity theorem \ref{t:1} relies
on this estimate.
\begin{proposition}\label{p}
Assume that $(H1)$, $(H2)$ and $(H3)$ hold. Let $ F $ be the map
defined in (\ref{eq:F2}) and assume that conditions (a)-(c) following
(\ref{eq:F2}) are true. Then we have the estimate
$$ 
\tb{u}_{\phi, \Omega_{2}} \lesssim \lambda^{1/2} \|Pu\|_{\phi,
    \Omega_{1}} + \tb{u}_{\phi, \Omega\setminus\Omega_{2}}
$$
where $\Omega$ is a neighbourood of $x_0$, $
\Omega_{2}\subset\!\subset \Omega_{1} \subset\!\subset \Omega 
$, $ \lambda \geq 1 $ suitably large
and for every $ u $ holomorphic in $\Omega $. 
\end{proposition}
The proof of the above result is split into several steps.   
We can decompose the set  $\Omega_2$ as follows  
$$
\Omega_2 = \Omega_{2,1} \cup \Omega_{2,2} \cup  \Omega_{2,3} 
$$ 
where 
$$
\Omega_{2,1}=\left \{ x\in  \Omega_2\, \Big|\, d_{1\phi} (x)\geq \frac 12 S
\lambda^{-1/2} \right \}, 
$$
$$
\Omega_{2,2}=\left \{ x\in  \Omega_2\, \Big|\, d_{1\phi} (x)\leq \frac 12 S
\lambda^{-1/2} \leq d_{2\phi} (x)\right \} 
$$
and
$$
\Omega_{2,3}=\left \{ x\in  \Omega_2\, \Big|\, (d_{1\phi}(x) \leq)
  \ d_{2\phi} (x)\leq \frac 12 S \lambda^{-1/2}\right \}. 
$$
We begin by localizing the problem in the ``elliptic'' region
$\Omega_{2,1}$.  
\begin{lemma}
\label{l.4.1} 
There exists a positive constant $ C $ such that
$$ 
\tb{u}_{\phi, \Omega_{2,1}} \lesssim  \left( \|Pu\|_{\phi,
    \Omega_{2,1}} +   \left( \frac 1S +e^{-\lambda /C}
\right)\tb{u}_{\phi, \Omega }\right) 
$$
with $ \lambda ,S\geq 1 $ suitably large and  for every $ u $ holomorphic in $
\Omega $. 
\end{lemma} 
{\sc Proof:}
We want to start by estimating 
\begin{multline*}
Mu(x)= \left ( \frac {\lambda}{2\pi }\right )^n \iint e^{i \lambda
  (x-x') \xi } \left (P(x,\xi , \lambda )-p_2\left ( x , \frac 2i
\partial_x \phi (x) \right )\right ) \\
u(x') \, dx'\wedge d\xi 
\end{multline*}
where the integration is performed along the contour
\begin{equation}\label{eq:contour} 
\xi = \frac{2}{i} \partial_{x}\phi(x) + i K (\overline{x-x'}),\qquad 
|x-x'| \leq r .
\end{equation}
Using the decomposition 
\begin{multline}
\label{eq:dec}
P(x,\xi , \lambda
)-p_2\left ( x , \frac 2i \partial_x \phi (x) \right )
\\
=
p_2(x,\xi )-p_2\left ( x , \frac 2i \partial_x \phi (x) \right ) +
P(x,\xi , \lambda
)-p_2( x , \xi )
\end{multline} 
We observe that, by Taylor formula, 
$$ 
\left | p_2(x,\xi )-p_2\left ( x , \frac 2i \partial_x \phi (x) \right
)\right | \lesssim  d_{1\phi }(x) K |x-x'| + (K |x-x'|)^2 .  
$$
Using the decomposition in equation (\ref{eq:dec}), we denote by $
M_{1} $ and $ M_{2} $ the corresponding pseudodifferential operators
in $ M $, so that
$$
Mu=M_1 u+M_2 u\, .
$$
In order to study the continuity of the operators $M_1$,
$M_2$ between the spaces $L^{2,2}_\phi ( \Omega_{2,1})$ and
$L^{2}_\phi ( \Omega_{2,1})$  it is
enough to estimate the corresponding reduced kernels. 

Let us preliminarily remark that, by Taylor expansion, 
\begin{multline*}
\frac {d_{1\phi }^{2}(x) + \lambda^{-1} }
{d_{1\phi }^{2}(x') + \lambda^{-1} }
\lesssim 1+ \frac {d_{1\phi } (x') |x-x'|+|x-x'|^2}{d_{1\phi }^{2}(x')
  + \lambda^{-1}}  
\\
\lesssim 1+\frac { |x-x'|}{d_{1\phi } (x')}+\left ( \frac {
  |x-x'|}{d_{1\phi } (x')}\right )^2.
\end{multline*}
Since $x,x'\in \Omega_{2,1}$ and $S\geq 1$, we find
\begin{multline*}
1+\frac { |x-x'|}{d_{1\phi } (x')}+\left ( \frac {
  |x-x'|}{d_{1\phi } (x')}\right )^2
\\
\lesssim \left (1+\frac{\lambda^{1/2}}{S} |x-x'|+ \frac{\lambda}{S^2}
|x-x'|^2\right )\lesssim
(1+\lambda^{1/2} |x-x'|)^2. 
\end{multline*}
Hence, using once more the fact that $x\in \Omega_{2,1}$ and
(\ref{eq:sal}), 
we find that the reduced kernel of $M_1$ can be estimated by 
\begin{multline*}
\frac{ d_{1\phi }(x) K |x-x'| + (K |x-x'|)^2}{ d_{1\phi }^{2}(x) +
  \lambda^{-1}} 
(1+\lambda^{1/2}|x-x'|)^2 \lambda^n e^{-\lambda K |x-x'|^2/C}
\\
\lesssim \left (K \frac {|x-x'|}{d_{1\phi }(x)}+ \left ( K  \frac
  {|x-x'|}{d_{1\phi }(x)}\right )^2\right ) (1+\lambda^{1/2}|x-x'|)^2
  \lambda^n e^{-\lambda K |x-x'|^2/C} 
\\
\lesssim 
\left ( \lambda^{1/2}S^{-1} K |x-x'|+ \lambda S^{-2} K^{2}
  |x-x'|^2\right )
\\
\times
(1+\lambda^{1/2}|x-x'|)^2 \lambda^n e^{-\lambda K |x-x'|^2/C}
\\
\lesssim 
\frac{K^{3/2}}{S} \lambda^{1/2} K^{1/2} |x-x'| 
(1+\lambda^{1/2}|x-x'|)^3 \lambda^n e^{-\lambda K |x-x'|^2/C}
\end{multline*}
i.e.
\begin{equation}
\label{eq:k1}
\| M_1 u \|_{\phi , \Omega_{2,1}} \lesssim  \frac{K^{3/2 - n}}{S}
\tb{u}_{\phi ,\Omega} \lesssim \frac{1}{S} \tb{u}_{\phi ,\Omega},
\end{equation}
since we may always have $ K \geq 1 $ and $ n \geq 3 $.

Let us estimate the reduced kernel of $ M_{2} $. 
For $\xi $ in the contour given by (\ref{eq:contour}),  
$$
|P(x,\xi ,\lambda )-p_2(x,\xi )|\lesssim \lambda^{-1} \left(
  d_{2\phi}(x) + K |x-x'|\right) \lesssim \lambda^{-1}
$$
Hence, using once more the fact that $x\in \Omega_{2,1}$ as in the
estimate of the reduced kernel of $M_1$, we find that the reduced
kernel of $M_2$ can be estimated by  
\begin{multline*}
\frac{\lambda^{-1}}{d_{1\phi}^{2}(x') + \lambda^{-1}} \lambda^n
  e^{-\lambda K |x-x'|^2/C} 
\\
\le \frac{\lambda^{-1}}{d_{1\phi }^{2}(x) + \lambda^{-1}} \lambda^n
(1+\lambda^{1/2}|x-x'|)^2  e^{-\lambda K |x-x'|^2/C}
\\
\lesssim  S^{-2} \lambda^n
(1+\lambda^{1/2}|x-x'|)^2  e^{-\lambda K |x-x'|^2/C}
\end{multline*}
so that
$$
\| M_2 u \|_{\phi , \Omega_{2,1}} \lesssim \frac{1}{S^{2}}
\tb{u}_{\phi ,\Omega} . 
$$ 
The above equation and (\ref{eq:k1}) yield 
\begin{equation}\label{eq:k}
\| M u \|_{\phi , \Omega_{2,1}} \lesssim \frac 1S \tb{u}_{\phi
  ,\Omega} . 
\end{equation} 
Set 
$$
Lu(x)=u(x)-\left ( \frac {\lambda}{2\pi}\right )^n\iint e^{i\lambda (x-x')\xi
} u(x') dx'd\xi  
$$
where the integral is once more performed along the contour in
(\ref{eq:contour}). Arguing as in \cite{Sj-Ast} we may show that 
$$
\| Lu\|_{\phi , \Omega_{2,1}} \lesssim e^{-\lambda /C} \tb{u}_{\phi ,
  \Omega} \, .
$$
Hence, by (\ref{eq:k}), we have  
\begin{multline}\label{eq:4.4}
\left \| Pu - p_2\left ( x , \frac 2i \partial_x \phi (x) \right ) u
\right \|_{\phi , \Omega_{2,1}} \le \| Mu \|_{\phi , \Omega_{2,1}}+\|
Lu \|_{\phi , \Omega_{2,1}}  
\\
\lesssim \left ( \frac 1S +
e^{-\lambda /C} \right )   \tb{u}_{\phi , \Omega }.  
\end{multline} 
Finally, we observe that, for $x\in \Omega_{2,1}$, 
$$
\left | p_2\left ( x , \frac
2i \partial_x \phi (x) \right )\right | \gtrsim d_{1\phi }^2 (x)
\gtrsim d_{1\phi }^2 (x) + \frac{S^2}{\lambda} \gtrsim d_{1\phi }^2
(x) + \frac{1}{\lambda} 
$$
hence
$$
 \left \|p_2\left ( x , \frac
2i \partial_x \phi (x) \right ) u 
\right \|_{\phi , \Omega_{2,1}} \gtrsim \tb{u}_{\phi , \Omega_{2,1} }.
$$
The above inequality and (\ref{eq:4.4}) yield the 
conclusion.\Qed

The following elementary covering result is needed 
to estimate the norm of $u$ in the set $\Omega_{2,3}=\{ x\in
\Omega_2\, |\,  d_{2\phi} (x)\leq (1/2)S 
\lambda^{-1/2}\}$.  

Given $x=(t,y,s)\in \C^n$ we define 
$$
V(x,S,\lambda )=B( (t,y), S\lambda^{-1/2})\times B(s, S\lambda^{-1/2}).
$$
\begin{lemma}
\label{l:cov}
There exist $ N > 0 $ and $N_0 < N$, $ N_{0} $ independent of
$S,\lambda$, such that, for every $S,\lambda\geq 1$ with 
$4S\lambda^{-1/2} <\min\{ r, \dist  (\Omega_2 , \complement \Omega_1)
\}$, we can find 
$$
x_1,\ldots ,x_N\in \pi_x (\Lambda_{\phi} \cap
\Sigma_2^{\C})\cap (\Omega_2 + V(0,S,\lambda ))
$$ 
such that 
\begin{itemize}
\item[(i)]
$\Omega_{2,3}\subset \cup_{j=1}^N V(x_j,  S , \lambda )$.
\item[(ii)]
$ N \lesssim \lambda^{n} $.
\item[(iii)] every point is contained in at most $N_0$ polydiscs
$ V(x_j,  2S , \lambda )$.   

\end{itemize}

\end{lemma} 
{\sc Proof:} 
We cut the $ x $-space into cubes with disjoint interior, $\{
Q_\alpha\}$, such that 
$$
\diam \pi_{t,y} Q_\alpha = S\lambda^{-1/2}\quad
\text{and}\quad 
\diam \pi_{s} Q_\alpha =S\lambda^{-1/2}. 
$$
In each so defined cube
intersecting the set $\pi_x (\Lambda_{\phi} \cap
\Sigma_2^{\C})\cap (\Omega_2 + V(0,S,\lambda ))$ we choose a
point $ x_{j} \in \pi_x (\Lambda_{\phi} \cap \Sigma_2^{\C})\cap
(\Omega_2 + V(0,S,\lambda))$, $ j = 1, \ldots, N $. 
Let $ Q_{j} $ denote the cube where the point $ x_{j} $ has been
picked, $ j=1, \ldots, N $.
Let $ x \in \Omega_{2,3} $. Then there is a point 
$ (\tilde{x},\tilde{\xi})=(0,0, \tilde{s} ,0,0, \tilde{\sigma} )\in
\Lambda_{\phi} \cap \Sigma_{2}^{\C} $  such that  
$$
d_{2\phi}(x) = \left| \left(x, \frac{2}{i} \partial_{x}\phi(x)\right) -
(\tilde{x},\tilde{\xi})\right| \leq \frac{S}{2\lambda^{1/2}}.
$$
Thus there exists a point $ x_{j}=(0,0,s_j)$, chosen above, such
that 
$$ 
|\tilde{s}-s_j|\leq  S\lambda^{-1/2} .
$$
We conclude that the polydiscs $
V(x_{j}, S, \lambda ) $, $ j =1, \ldots, N $, are a covering of $
\Omega_{2,3} $. 

Now $ \cup_{j=1}^{N}V(x_{j}, S ,\lambda ) \subset \Omega_2 +
V(0, 2S,\lambda ) $ and $ Q_{j} \subset V (x_{j}, S,\lambda )
$. Since, by assumption, the volume of $ \Omega_2 +
V(0,2S,\lambda ) $ is bounded by a constant independent of $ S $
and $ \lambda $, because 
$$
\vol \left(\Omega_2 + V(0,2S,\lambda ) \right) \geq \vol \left(
  \cup_{j=1}^{N} Q_{j}\right) \gtrsim N \lambda^{-n},
$$
we conclude (ii).

Moreover, if a point belongs to $ V(x_j, 2S,\lambda )\cap  V(x_k,
 2S,\lambda )$ then $x_k\in  V(x_j, 4S, \lambda )$. 
Hence, slightly  enlarging the polydisc $V(x_j, 4S,\lambda )$, we can
suppose that the whole cube $ Q_{k} $, containing
$x_k$, is a subset of that polydisc. 
Since the so enlarged polydisc may contain at most a finite number of
cubes $ Q_{\ell} $ and both the volume of the enlarged polydisc and that
of $ Q_{\ell} $ is $ \mathcal{O}(S^{2n} \lambda^{-n}) $, we obtain
the third item in the statement.
\Qed 
\begin{remark}
We observe that Condition (iii) above implies the following
equivalence of norms (with constants independent of $\lambda$ and $S$)  
\begin{equation}
\label{eq:eqn}
\begin{array}{l}
\| \cdot \|_{\phi , \cup_{j=1}^N V(x_j,2S,\lambda )}\sim \sum_{j=1}^N
\| \cdot \|_{\phi , V(x_j,2S,\lambda )} 
\\[12pt]
\tb{ \cdot}_{\phi , \cup_{j=1}^N V(x_j,2S,\lambda )}\sim \sum_{j=1}^N
\tb{\cdot }_{\phi , V(x_j,2S,\lambda )}.  
\end{array} 
\end{equation} 
\end{remark}
Next we prove an a priori estimate in $ \Omega_{2,3} $.
\begin{lemma}
\label{l.4.2} 
There exists a positive function  $ \tilde{\gamma}(S,K,\lambda) $ such that
$$ 
\tb{u}_{\phi, \Omega_{2,3}} \lesssim (1 + \mathcal{O}_{S}(\epsilon))
\left( \lambda^{1/2}\|Pu\|_{\phi,
    \Omega_{2,3}} +  \tilde{\gamma}(S,K,\lambda )\tb{u}_{\phi, \Omega }\right) 
$$
for every $ u $ holomorphic in $\Omega $. 
$$
\tilde{\gamma}(S,K,\lambda )=   e^{-\lambda /C} +\frac{S^3} {K^n}
+ K^{3/2-n}+ \lambda^{-1/2} S^{n_{t} +4} +e^{-S/C}. 
$$
\end{lemma}
{\sc Proof:} Let $z\in \Omega_2$, $(z, \zeta )=(0,0,z_s,
0,0,\zeta_\sigma )\in  \Lambda_{\phi }
\cap \Sigma_2^{\C}$ . We have 
\begin{multline*}
| P(x,\xi ,\lambda )- (P)_{(z,\zeta )} (x,\xi ,\lambda )| 
\\[5pt]
\lesssim 
\Big(d_{1\phi}^2(x) + (K |x-x'|)^{2}\Big) \Big[ d_{1\phi}(x) +d_{2\phi}(x)
  +|s-z_s| + K |x - x'| \Big]
\\[5pt]
+\lambda^{-1} \Big(d_{2\phi}(x) +K |x - x'|\Big) \Big[d_{1\phi}(x) +d_{2\phi}(x)
  +|s-z_s| + K |x - x'|  \Big]
\\[5pt]
+ \mathcal{O}\left(\lambda^{-2}\right),
\end{multline*} 
where $ \mathcal{O}\left(\lambda^{-2}\right) $ is uniform with respect
to $ S $ and $ K $.
We recall that in the above estimate we used that, for every $\xi$ in
the contour given in (\ref{eq:contour}), we have
$$
|\xi -\zeta |\lesssim |x-z|+ K |x-x'|, 
$$
uniformly with respect to $ S $.
We realize $(P)_{(z,\zeta )}$ as a pseudodifferential operator,
integrating along the same contour used for the realization of $P$. 
In order to study the continuity of the operator 
$$ 
P-(P)_{(z,\zeta )}: L^{2,2}_\phi (\Omega )\longrightarrow L^{2}_\phi (
V(z,2S,\lambda )) 
$$
we remark that the corresponding reduced kernel can be estimated 
(modulo constants) by 
\begin{multline}
\label{eq:redk}
\frac {\lambda^n e^{-\lambda K |x-x'|^2/C}}{     
d_{1\phi}^2 (x')+ \lambda^{-1}} \times
\\[5pt]
\left\{
\Big(d_{1\phi}^2(x) + (K |x-x'|)^{2}\Big) \Big[ d_{1\phi}(x) +d_{2\phi}(x)
  +|s-z_s| + K |x - x'| \Big] \right.
\\[5pt]
+\lambda^{-1} \Big(d_{2\phi}(x) +K |x - x'|\Big) \times
\\[5pt]
\left.
\Big[d_{1\phi}(x) +d_{2\phi}(x)
  +|s-z_s| + K |x - x'|  \Big]
+ \mathcal{O}\left(\lambda^{-2}\right)
\right\}
\\[5pt]
\lesssim
\frac {\lambda^n e^{-\lambda K |x-x'|^2/C}}{     
d_{1\phi}^2 (x')+ \lambda^{-1}} \left\{ \frac{S^{3}}{\lambda^{3/2}} +
  (K |x - x'|)^{3} + \frac{1}{\lambda} \left( \frac{S^{2}}{\lambda} +
    K^{2} |x-x'|^{2}\right)
\right\}.
\end{multline}
Hence
\begin{multline}
\label{eq:eee}
\lambda^{1/2} \| (P-(P)_{(z,\zeta )}) u \|_{\phi,V(z, 2S, \lambda)}
\\
\lesssim \left\{ \frac{S^{3}}{K^{n}} + \mathcal{O}\left(K^{\frac{3}{2}
    - n}\right) +
\mathcal{O}\left(\frac{S^{2}}{\lambda^{1/2}}\right)\right\} 
\tb{u}_{\phi,V(z, 2S, \lambda)}.
\end{multline}
Let now $x_1,\cdots , x_N\in  \Pi_x (\Lambda_{\phi} \cap
\Sigma_2^{\C})$ be the points given in Lemma \ref{l:cov}. 
Using the estimate (\ref{eq:eee}) and the equivalence of the norms in
(\ref{eq:eqn}), we arrive at the following inequality
\begin{multline}
\label{eq:normg} 
\lambda^{1/2} \sum_{j=1}^N \| ((P)_{(x_j,\xi_j)}-P )u\|_{\phi,V(x_j,
    2S,\lambda)} 
\\
\lesssim  N_{0} \left\{ \frac{S^{3}}{K^{n}} +
    \mathcal{O}\left(K^{\frac{3}{2} - n}\right) +
\mathcal{O}\left(\frac{S^{2}}{\lambda^{1/2}}\right)\right\} 
\tb{u}_{\phi,\cup_{j=1}^{N} V(x_{j}, 2S, \lambda)}.
\end{multline}
We use the notation $P_{(x_j,\xi_j )}u=P_{(x_j,\xi_j )}(x,\tilde{D},\lambda
)u$, $j=1,\cdots ,N$, for the action of the operator $(P)_{(x_j,\xi_j
  )}$ as a differential operator on the function $u$. 
We recall that, in \cite{Sj-Ast}, the following estimate 
is proved 
\begin{equation}
\label{eq:locdif}
\| (P)_{(x_j,\xi_j )}u - P_{(x_j,\xi_j )}u\|_{\phi,V(x_j,
  2S,\lambda)}\lesssim e^{-\lambda /C} \tb{u}_{\phi , \Omega }.  
\end{equation}
The Estimates (\ref{eq:normg}) and (\ref{eq:locdif}) yield  
\begin{multline}
\label{eq:4.21}
\lambda^{1/2} \sum_{j=1}^N \| P_{(x_j,\xi_j)}u\|_{\phi , V(x_j, 2S,\lambda)} 
\\
\leq \lambda^{1/2} \sum_{j=1}^N \| ( P_{(x_j,\xi_j)}-
(P)_{(x_j,\xi_j)} ) u\|_{\phi , V(x_j, 2S,\lambda)}
\\
+\lambda^{1/2} \sum_{j=1}^N \Big[\| ((P)_{(x_j,\xi_j)}-P )u\|_{\phi , V(x_j,
    2S,\lambda)}+ \| P u\|_{\phi , V(x_j, 2S,\lambda)} \Big]
\\
\lesssim  \left(\lambda^{1/2} N e^{-\lambda/C} +  N_{0} \left\{
    \frac{S^{3}}{K^{n}} + 
    \mathcal{O}\left(K^{\frac{3}{2} - n}\right) +
\mathcal{O}\left(\frac{S^{2}}{\lambda^{1/2}}\right)\right\} \right)
\tb{u}_{\phi , \Omega }
\\
+\lambda^{1/2} \| P u\|_{\phi ,\Omega_{1} }.
\end{multline}

Now, we want to show that for every $S\geq 1$ and for $\lambda$ large
if $F$ is close enough to the identity map then  
\begin{multline}\label{eq:lemma}
\tb{u}_{\phi , V(x_j,
  S,\lambda)}\lesssim \lambda^{1/2} \| P_{(x_j,\xi_j )}u
  \|_{\phi , V(x_j, 
  2S,\lambda )} 
\\
+ \gamma (S,\lambda )\tb{u}_{\phi , V(x_j,
  2S,\lambda)} 
\end{multline}
for every $j=1,\ldots ,N$ and for every $u$ holomorphic in $V(x_j,
  2S,\lambda)$.  
Here  
$$ 
\gamma(S, \lambda)= e^{-S/C} +\lambda^{-1/2} S^{n_{t} +4}. 
$$
Essentially, Estimate (\ref{eq:lemma}) reduces to the estimate in
Proposition \ref{p:2} for the localized operator.

Indeed by (\ref{eq:4.5}), 
\begin{multline*}
\tb{u}_{\phi, V(x_{j}, S, \lambda)} = \lambda^{-\frac{n}{2}+1}
e^{-\lambda\phi(x_{j})}) (1 + \mathcal{O}_{S}(\epsilon)) \tb{v}_{\phi_{0},
B_{(t, y)}(0, S) \times B_{s}(0, S)} \\
\lesssim \lambda^{-\frac{n}{2}+1}
e^{-\lambda\phi(x_{j})}) (1 + \mathcal{O}_{S}(\epsilon)) \times
\\
\Big(
  \lambda^{1/2} \| [P]_{(x_{j}, \xi_{j})}v\|_{\phi_{0}, B_{(t,
      y)}(0,2S)\times B_{s}(0, 2S)} 
\\
+ \gamma(S, \lambda) \tb{v}_{\phi_{0}, B_{(t,
      y)}(0,2S)\times B_{s}(0, 2S)}  \Big)
\\
\lesssim 
(1 + \mathcal{O}_{S}(\epsilon))
\left(
  \lambda^{1/2} \| P_{(x_{j}, \xi_{j})}u\|_{\phi_{0}, V(x_{j}, 2S,
    \lambda)} + \gamma(S, \lambda) \tb{u}_{\phi_{0}, V(x_{j}, 2S, \lambda)}  \right).
\end{multline*}
In the above inequalities we used Proposition \ref{p:2} and formula
(\ref{eq:4.7}). Moreover here $ v $ has been defined in (\ref{eq:U}).

Now, we have 
\begin{multline*}
\sum_{j=1}^N \tb{u}_{\phi , V(x_j,S, \lambda)}
\\
\lesssim (1 + \mathcal{O}_{S}(\epsilon))
\sum_{j=1}^N  \Big( \lambda^{1/2} \| P_{(x_j,\xi_j )}u
  \|_{\phi , V(x_j,2S,\lambda)} + \gamma (S,\lambda )\tb{u}_{\phi , V(x_j,
  2S,\lambda)}  \Big)
\\
\lesssim (1 + \mathcal{O}_{S}(\epsilon))
\Big( \lambda^{1/2} \| Pu\|_{\phi , \Omega_{1}} + G(K, S, \lambda)
 \tb{u}_{\phi,\Omega }\Big).
\end{multline*}
Here
\begin{multline}
\label{eq:G}
G(K, S, \lambda) =  N_{0} \left\{
    \frac{S^{3}}{K^{n}} + 
    \mathcal{O}\left(K^{\frac{3}{2} - n}\right) +
\mathcal{O}\left(\frac{S^{2}}{\lambda^{1/2}}\right)\right\}
\\
+ \lambda^{1/2} N e^{-\lambda/C}  + e^{-S/C} + \lambda^{-1/2} S^{n_{t}+4}
\end{multline}
is the quantity defined in (\ref{eq:4.21}) and $ \gamma(S, \lambda) $
defined after (\ref{eq:lemma}).
This completes the proof of the lemma. 
\Qed 

The next lemma takes care of the microlocal region $ \Omega_{2,2} $.
\begin{lemma}
\label{l.4.3} 
There exists a positive function  $ \tilde{\gamma}(S,\lambda) $ such that
$$ 
\tb{u}_{\phi, \Omega_{2,2}} \lesssim \lambda^{1/2}\|Pu\|_{\phi,
    \Omega_{2,2}} +  \tilde{\gamma}(S, K, \lambda ) \tb{u}_{\phi, \Omega } 
$$
for every $ u $ holomorphic in $\Omega $. 
\end{lemma}
{\sc Proof:}
The proof is done following the same ideas of the proof of Lemma
\ref{l.4.2}; we just sketch it out emphasizing the main
differences. As above we cover $ \Omega_{2,2} $ with polydiscs
centered at points of $
\big(\Sigma_{1}^{\C}\setminus\Sigma_{2}^{\C}\big) \cap \Lambda_{\phi}
$. The basic ingredient in the proof is an a priori estimate for the
operator localized at each polydisc center. Once this estimate is
obtained we use the perturbation argument in the proof of Lemma
\ref{l.4.2} to get rid of the error terms using the large parameters $
S $, $ K $ and $ \lambda $. We would like to stress the fact that the
size of $ S $ and $ K $ at this stage, as well as at the previous
stage, depends only on the problem's data.

The localized operator is given by (\ref{eq:loc22}). For this operator
an approximate parametrix can be constructed along the same lines of
Section \ref{s:3}. The only difference in the present case is that the
``lower order terms'' are elliptic. On the other hand, due to the fact
that we are in the region $ \Omega_{2,2} $, the lower order
term can be estimated from below by $ C S \lambda^{-1/2} $ and this
is responsible of a factor $ \lambda^{1/2} $ in the ensuing a priori
estimate. 

\Qed

We are now ready to prove the main a priori estimate.

{\sc Proof of the Proposition \ref{p}:} Using the Lemmas \ref{l.4.1},
\ref{l.4.2} and \ref{l.4.3} we have 
$$
\tb{u}_{\phi, \Omega_{2}}\leq C (1 + \mathcal{O}_{S}(\epsilon)) 
\Big( \lambda^{1/2}\|Pu\|_{\phi,
    \Omega_{1}} +  \tilde{\gamma}(S, K, \lambda )\tb{u}_{\phi, \Omega
  } \Big),
$$
for a suitable positive constant $ C $ independent of $ S $, $ K $, $
\lambda $ and $ \epsilon $. Furthermore 
$$
\tilde{\gamma}(S,K,\lambda )=   e^{-\lambda /C} +\frac{S^3} {K^n}
+ K^{3/2-n}+ \lambda^{-1/2} S^{n_{t} +4} +e^{-S/C} + \frac{1}{S}. 
$$
We recall that the parameter $ S $ has to be chosen large but
depending only on the given operator. Hence the quantity 
$ C (1+\mathcal{O}_{S}(\epsilon)) \tilde{\gamma}(S, K, \lambda) $ can
be made smaller than $ 1/2 $ choosing $ S $ large and $ \epsilon $
small depending only on the operator, $ K $ large
depending on $ S $ and $ \lambda $ suitably large.

Hence,    
$$
\tb{u}_{\phi, \Omega_{2}}\lesssim \lambda^{1/2}\|Pu\|_{\phi,
    \Omega_{1}} +  \tb{u}_{\phi,
  \Omega\setminus \Omega_2  }
$$
and the proof is completed. 
\Qed

\section{The construction of the phase $\phi$}
\renewcommand{\theequation}{\thesection.\arabic{equation}}
\setcounter{equation}{0}
\setcounter{theorem}{0}
\setcounter{proposition}{0}  
\setcounter{lemma}{0}
\setcounter{corollary}{0} 
\setcounter{definition}{0}
\subsection{Remarks on the Hamilton-Jacobi equation}

Let $ (x_{0}, \xi_{0}) \in \Sigma_{2} $, $ W $ be an open neighborhood
of $ (x_{0}, \xi_{0}) $ in $ \C^{2n} $ and $ \Omega \subset W $.

The weight function $ \phi $ is constructed by solving for small
values of the time variable $ t $ a Hamilton-Jacobi equation.

Let $ r \colon W \rightarrow \C $ be a $ C^{\infty} $
function. Consider
$$ 
\begin{cases}
\displaystyle{   \frac{\partial \phi}{\partial t}(t, x)} &
    \displaystyle{= (\re r) \left(x, 
    \frac{2}{i}\partial_{x} \phi(t, x)\right)}
\\[16pt]
\phi(0, x) & = \phi_{0}(x),
\end{cases}
$$
for $ 0 \leq t \leq \epsilon_{0} $.

The solution of the above problem is constructed using the standard
Hamilton-Jacobi theory with respect to the symplectic form 
$$ 
\im\sigma = \im \left(d\xi \wedge dx\right).
$$
Actually, setting $ \phi_{t}(x) = \phi(t, x) $, we have
$$ 
\Lambda_{\phi_{t}} = \exp\left( t H_{\scriptscriptstyle \re
    r}^{\scriptscriptstyle\im \sigma} \right) \Lambda_{\phi_{0}}. 
$$
If $ r $ is a holomorphic function on $ W $ we have
$$ 
H_{\scriptscriptstyle \re r}^{\scriptscriptstyle\im \sigma} = \widehat{H_{ir}},
$$
where $H_{ir} $ is the usual complex standard Hamilton field of $ ir $
and $ \widehat{H_{ir}} $ denotes the real part of $ H_{ir} $, i.e. the
real field that gives the same result as $ H_{ir} $ when acting on
holomorphic functions.

\begin{remark}
If $ r $ is holomorphic in $ W $ and real valued on $
\Lambda_{\phi_{0}} $ the solution of the above Hamilton-Jacobi problem
is obtained as the restriction to the positive $ t $-axis of the
solution of the complex equation
$$ 
\begin{cases}
\partial_{t} \psi (t,x)& = r\left( x, \frac{2}{i}\partial_{x}\psi(t,
  x)\right) \\[16pt]
\psi(0, x) & = \phi_{0}(x),
\end{cases}
$$
for $ |t| < \epsilon_{0} $.
\end{remark}

\subsection{Contruction of the function $ r $}

Since $ \R^{2n} $ and $ \Lambda_{\phi_{0}} $ are isomorphic it is
easier to contruct the function $ r $ in $ \R^{2n} $ near the point 
$$
(x_{0}, \xi_{0}) = (0, 0, s_{0}; 0, 0, \sigma_{0}) \in
\tilde{\Sigma}_{2},
$$
where $ \tilde{\Sigma}_{2} $ is the (real) characteristic manifold of
the first eigenvalue of $ \tilde{P} $.

We want
$$ 
H_{r}(\rho_{j}) \in T\tilde{\Sigma}_{j}\, , \qquad \text{if} \ \rho_{j} \in
\tilde{\Sigma}_{j}. 
$$

Let us choose
$$
r(x, \xi) = (s - s_{0})^{2} + (\sigma - \sigma_{0})^{2} 
+ C(|t|^{2} + |\tau|^{2} + y^{2} + \eta^{2}),
$$
where $ C $ is a positive constant that makes $ r $ as positive as we
desire outside $ \tilde{\Sigma}_{2} $.

We have that in the real domain $ \R^{2n} $ 
$$ 
r (x, \xi) \sim |x - x_{0}|^{2} + |\xi - \xi_{0}|^{2}.
$$
Then, on $ \Lambda_{\phi_{0}} $,
$$ 
r\left(x, \frac{2}{i}\partial_{x}\phi_{0}\right) \sim |x - x_{0}|^{2},
$$
for every $ x \in \pi_{x}(W \cap \Lambda_{\phi_{0}}) $.

\medskip

\section{Proof of Theorem \ref{t:1}}
\renewcommand{\theequation}{\thesection.\arabic{equation}}
\setcounter{equation}{0}
\setcounter{theorem}{0}
\setcounter{proposition}{0}  
\setcounter{lemma}{0}
\setcounter{corollary}{0} 
\setcounter{definition}{0}

We want to show that if $(x_0,\xi_0)\notin WF_a(\tilde{P}u)$ then 
$(x_0,\xi_0)\notin WF_a(u)$. 
We recall the \textit{a priori} estimate obtained:
$$ 
\tb{u}_{\phi, \Omega_{2}} \leq C \left( \lambda^{1/2} \|Pu\|_{\phi,
    \Omega_{1}} + \tb{u}_{\phi, \Omega\setminus\Omega_{2}}\right)
$$
where $ \Omega_{2}\subset\!\subset \Omega_{1} \subset\!\subset \Omega
\subset W$ and $ x_{0} \in \Omega_{2} $ and we write $u$ and $x_0$
instead of $Tu$ 
and $x_0-i\xi_0$ respectively.

Since $ u $ is a tempered distribution before the FBI transform, we
have 
$$ 
\| u \|_{\phi_{0},\Omega} \leq C \lambda^{N_{0}},
$$
for a certain $ N_{0} \in \N $.

Since $ \tilde{P}u $ is real analytic at the real point $ (x_{0},
\xi_{0}) $ 
before the FBI transform, we have
$$ 
\| Pu \|_{\phi_{0}, \Omega_{3}} \leq C_{1} e^{-\lambda/C_{1}},
$$
for a positive constant $ C_{1} $; here $ \Omega_{3} $ is a suitable
neighborhood of $ x_{0} $. Recalling that 
$$ 
\phi_{t}(x) - \phi_{0}(x) \sim t |x - x_{0}|^{2},
$$
we obtain that
$$ 
\| Pu \|_{\phi_{t}, \Omega_{1}} \leq \tilde{C} e^{-\lambda/\tilde{C}}, 
$$
for a positive constant $ \tilde{C} $.

Decompose $ \Omega \setminus \Omega_{2} = K_{1} \cup K_{2} $, where 
$$ 
K_{2} \cap \pi_{x}(\Sigma_{2}) = \varnothing,
$$
while
$$ 
r_{|_{K_{1}}} \geq \alpha > 0.
$$

Since
$$ 
{\phi_{t}}_{|_{K_{1}}} \geq \phi_{0} + \alpha_{1} t, \qquad \alpha_{1} > 0,
$$
we have
$$ 
\tb{u}_{\phi_{t},K_{1}} \leq C e^{-\lambda/C_{t}'}
\|u\|_{\phi_{0},\Omega_{1}} \leq C_{t} e^{-\lambda/C_{t}}
$$

Since $ \tilde{P} $ is analytic hypoelliptic far from $
\tilde{\Sigma}_{2} $, by the Tartakoff-Treves theorem, we have
$$ 
\tb{u}_{\phi_{0},K_{2}} \leq C_{2} e^{-\lambda/C_{2}}, \qquad C_{2} > 0.
$$

Arguing as above we get
$$ 
\tb{u}_{\phi_{t},K_{2}} \leq C_{2} e^{-\lambda/C_{2}}, \qquad C_{2} > 0.
$$

Hence the \textit{a priori} estimate implies that
$$ 
\tb{u}_{\phi_{t},\Omega_{2}} \leq C e^{-\lambda/C}.
$$
Let now $ \Omega_{4} $ be a sufficiently small neighborhood of $ x_{0}
$ such that
$$ 
\phi_{t}(x) < \phi_{0}(x) + \frac{1}{3C_{4}}
$$
on $ \Omega_{4} $. Then
$$ 
\|u \|_{\phi_{0},\Omega_{4}} \leq C e^{-\lambda/C},
$$
which means that $ u $ is real analytic at $ (x_{0}, \xi_{0}) $ before
the FBI transform. This proves the theorem.

\section{Some related model operators}
\renewcommand{\theequation}{\thesection.\arabic{equation}}
\setcounter{equation}{0}
\setcounter{theorem}{0}
\setcounter{proposition}{0}  
\setcounter{lemma}{0}
\setcounter{corollary}{0} 
\setcounter{definition}{0}

We briefly discuss in this section a case related to what we study
in the paper. The lowest eigenvalue in this case is identically zero
on a ``half fiber'' over a characteristic point.

Let us consider the operator $ \Box_{b} $ on functions for the
Heisenberg vector fields (strongly pseudo convex case). We use the
following notation:
$$ 
w = x+iy\in \C^n, \quad W_{j} = \partial_{w_{j}} + i \bar{w}_{j}
\frac{\partial}{\partial t}.
$$
$$ 
\Box_{b} = -\frac{1}{2} \sum_{j=1}^{n} \left( W_{j} \bar{W}_{j} +
  \bar{W}_{j} W_{j}\right) + i n \frac{\partial}{\partial t}.
$$
We have that $ \Char \Box_{b} = \{(0, 0, t; 0, 0, \tau)\}$ ,
$\tr^{+}\Box_{b} = n |\tau|$  
and
$$ 
(\Box_{b})^{s} + \tr^{+}\Box_{b} = n (|\tau| - \tau)\, .
$$
\begin{proposition}
Let $u\in \mathcal{D}'$ be a solution of the equation
\begin{equation}\label{eq:eq} 
\Box_{b}u = f
\end{equation}
where $ f \in C^{\omega}(U) $, $ U $ is an open set containing the point $
(0, 0, t_{0}) $. Assume that, there exists $\tau_{0} > 0$ such that  
\begin{equation}\label{eq:h} 
(0, 0, t_{0}; 0, 0, \tau_{0}) \notin WF_{a}(u)\, .
\end{equation}
Then $ u $ is real analytic at $ (0, 0, t_{0}) $.
\end{proposition}
In other words, to get analytic regularity of a solution of Equation
(\ref{eq:eq}), we need to 
assume that there are no analytic singularities of $u$  
in the region where $(\Box_{b})^{s} + \tr^{+}\Box_{b}$
identically vanishes.
\begin{remark}
In particular the same result of the above proposition holds for the
operator 
$$ 
P_{-} = D_{x}^{2} + x^{2} D_{t}^{2} - D_{t}.
$$
The result is optimal. In fact we have solutions of the equation $
P_{-}u = 0 $ with $ u \in G^{s} $ and no better for every $ s > 1 $.
It is enough to consider the function 
$$ 
u(x, t) = \int_{0}^{+\infty} e^{-\rho \frac{x^{2}}{2} + i\rho t -
  \rho^{1/s}} d\rho.
$$
and observe that 
$$ 
D_{t}^{k} u (0, 0) = \int_{0}^{+\infty} \rho^{k}
  e^{-\rho^{1/s}} d\rho \sim k!^{s} 
$$
\end{remark}
{\sc Proof:} We want to show that 
$$
( 0, 0, t_{0}; 0, 0, \tau ) \notin WF_{a}(u) \ \text{ when }\
\tau \not= 0\, . 
$$
Since $WF_{a}(u)$ is a conic set, by Assumption (\ref{eq:h}), we get 
$$
( 0, 0, t_{0}; 0, 0, \tau ) \notin WF_{a}(u) \ \text{ when }\
\tau >  0\, . 
$$ 
Moreover, in the set $\{ (0, t_{0}; 0, 0, \tau )\, |\, \tau <0 \}$
we have 
$$
(\Box_{b})^{s} + \tr^{+}\Box_{b}>0 
$$
hence the conclusion follows by the Tartakoff-Treves theorem.
\Qed

\end{document}